\documentclass{article}
\usepackage{amsmath}
\usepackage{amssymb}
\usepackage{amsthm}
\usepackage{graphicx}
\usepackage{cite}
\usepackage[arrow,curve,matrix,arc]{xy}
\DeclareFontFamily{U}{rsfs}{} \DeclareFontShape{U}{rsfs}{n}{it}{<->
rsfs10}{} \DeclareSymbolFont{mscr}{U}{rsfs}{n}{it}
\DeclareSymbolFontAlphabet{\scr}{mscr}
\def\mathscr{\scr}
\begin{document}
\def\e#1\e{\begin{equation}#1\end{equation}}
\def\ea#1\ea{\begin{align}#1\end{align}}
\def\eq#1{{\rm(\ref{#1})}}
\theoremstyle{plain}
\newtheorem{thm}{Theorem}[section]
\newtheorem{prop}[thm]{Proposition}
\newtheorem{cor}[thm]{Corollary}
\newtheorem{conj}[thm]{Conjecture}
\theoremstyle{definition}
\newtheorem{dfn}[thm]{Definition}
\newtheorem{ex}[thm]{Example}
\def\dim{\mathop{\rm dim}\nolimits}
\def\End{\mathop{\rm End}\nolimits}
\def\rank{\mathop{\rm rank}}
\def\Hilb{\mathop{\rm Hilb}\nolimits}
\def\coh{\mathop{\rm coh}\nolimits}
\def\cs{{\rm cs}}
\def\id{\mathop{\rm id}\nolimits}
\def\bdim{{\mathbin{\bf dim}\kern.1em}}
\def\Crit{\mathop{\rm Crit}}
\def\CF{\mathop{\rm CF}\nolimits}
\def\GL{\mathop{\rm GL}}
\def\U{\mathop{\ts\rm U}}
\def\Gr{\mathop{\rm Gr}}
\def\Mo{\mathop{\text{M\"o}}}
\def\Tr{\mathop{\rm Tr}}
\def\ch{\mathop{\rm ch}\nolimits}
\def\num{{\rm num}}
\def\vir{{\rm vir}}
\def\stp{{\rm stp}}
\def\fr{{\rm fr}}
\def\stf{{\rm stf}}
\def\stk{{\rm stk}}
\def\ss{{\rm ss}}
\def\st{{\rm st}}
\def\vi{{\rm vi}}
\def\an{{\rm an}}
\def\Aut{\mathop{\rm Aut}}
\def\Hom{\mathop{\rm Hom}\nolimits}
\def\Iso{\mathop{\rm Iso}\nolimits}
\def\Ext{\mathop{\rm Ext}\nolimits}
\def\Spec{\mathop{\rm Spec}\nolimits}
\def\LCF{\mathop{\rm LCF}\nolimits}
\def\SF{\mathop{\rm SF}\nolimits}
\def\SFa{\mathop{\rm SF}\nolimits_{\rm al}}
\def\SFai{\mathop{\rm SF}\nolimits_{\rm al}^{\rm ind}}
\def\oSF{\mathop{\bar{\rm SF}}\nolimits}
\def\oSFa{\mathop{\bar{\rm SF}}\nolimits_{\rm al}}
\def\oSFai{{\ts\bar{\rm SF}{}_{\rm al}^{\rm ind}}}
\def\modCQ{\text{\rm mod-$\C Q$}}
\def\modCQI{\text{\rm mod-$\C Q/I$}}
\def\bs{\boldsymbol}
\def\ge{\geqslant}
\def\le{\leqslant\nobreak}
\def\sA{{\mathbin{\mathscr A}}}
\def\sB{{\mathbin{\mathscr B}}}
\def\sG{{\mathbin{\mathscr G}}}
\def\bG{{\mathbin{\mathbb G}}}
\def\A{{\mathbin{\mathcal A}}}
\def\B{{\mathbin{\mathcal B}}}
\def\fE{{\mathbin{\mathfrak E}}}
\def\fF{{\mathbin{\mathfrak F}}}
\def\fG{{\mathbin{\mathfrak G}}}
\def\fH{{\mathbin{\mathfrak H}}}
\def\R{{\mathbin{\mathbb R}}}
\def\Z{{\mathbin{\mathbb Z}}}
\def\cQ{{\mathbin{\mathcal Q}}}
\def\Q{{\mathbin{\mathbb Q}}}
\def\N{{\mathbin{\mathbb N}}}
\def\C{{\mathbin{\mathbb C}}}
\def\K{{\mathbin{\mathbb K}}}
\def\CP{{\mathbin{\mathbb{CP}}}}
\def\M{{\mathbin{\mathcal M}}}
\def\fM{{\mathbin{\mathfrak M}}}
\def\O{{\mathcal O}}
\def\fR{{\mathbin{\mathfrak R}}}
\def\fS{{\mathbin{\mathfrak S}}}
\def\fExact{\mathop{\mathfrak{Exact}\kern .05em}\nolimits}
\def\fVect{{\mathbin{\mathfrak{Vect}}}}
\def\fW{{\mathbin{\mathfrak W}}}
\def\al{\alpha}
\def\be{\beta}
\def\ga{\gamma}
\def\de{\delta}
\def\bde{\bar\delta}
\def\bdss{\bar\delta_{\rm ss}}
\def\io{\iota}
\def\ep{\epsilon}
\def\bep{\bar\epsilon}
\def\la{\lambda}
\def\ka{\kappa}
\def\th{\theta}
\def\ze{\zeta}
\def\up{\upsilon}
\def\vp{\varphi}
\def\si{\sigma}
\def\om{\omega}
\def\La{\Lambda}
\def\Om{\Omega}
\def\Ga{\Gamma}
\def\pd{\partial}
\def\db{{\bar\partial}}
\def\ts{\textstyle}
\def\sst{\scriptscriptstyle}
\def\w{\wedge}
\def\sm{\setminus}
\def\bu{\bullet}
\def\op{\oplus}
\def\ot{\otimes}
\def\ov{\overline}
\def\bigop{\bigoplus}
\def\bigot{\bigotimes}
\def\iy{\infty}
\def\es{\emptyset}
\def\ra{\rightarrow}
\def\ab{\allowbreak}
\def\longra{\longrightarrow}
\def\hookra{\hookrightarrow}
\def\t{\times}
\def\ci{\circ}
\def\ti{\tilde}
\def\d{{\rm d}}
\def\md#1{\vert #1 \vert}
\title{Generalized Donaldson--Thomas invariants}
\author{Dominic Joyce}
\date{}
\maketitle
\begin{abstract}
This is a survey of the book \cite{JoSo} with Yinan Song.
Donaldson--Thomas invariants $DT^\al(\tau)\in\Z$ `count'
$\tau$-(semi)stable coherent sheaves with Chern character $\al$ on a
Calabi--Yau 3-fold $X$. They are unchanged under deformations of
$X$. The conventional definition works only for classes $\al$ with
no strictly $\tau$-semistable sheaves. Behrend showed that
$DT^\al(\tau)$ can be written as a weighted Euler characteristic
$\chi\bigl(\M_\st^\al(\tau), \nu_{\M_\st^\al(\tau)}\bigr)$ of the
stable moduli scheme $\M_\st^\al(\tau)$ by a constructible function
$\nu_{\M_\st^\al(\tau)}$ we call the `Behrend function'.

We discuss {\it generalized Donaldson--Thomas invariants\/}
$\bar{DT}{}^\al(\tau)\in\Q$. These are defined for all classes
$\al$, and are equal to $DT^\al(\tau)$ when it is defined. They are
unchanged under deformations of $X$, and transform according to a
known wall-crossing formula under change of stability condition
$\tau$. We conjecture that they can be written in terms of integral
{\it BPS invariants\/} $\hat{DT}{}^\al(\tau)\in\Z$ when the
stability condition $\tau$ is `generic'.

We extend the theory to abelian categories $\modCQI$ of
representations of a quiver $Q$ with relations $I$ coming from a
superpotential $W$ on $Q$, and connect our ideas with Szendr\H oi's
noncommutative Donaldson--Thomas invariants, and work by Reineke and
others on invariants counting quiver representations. The book
\cite{JoSo} has significant overlap with a recent, independent paper
of Kontsevich and Soibelman~\cite{KoSo1}.
\end{abstract}

\section{Introduction}
\label{ig1}

This is a survey of the book \cite{JoSo} by the author and Yinan
Song. Let $X$ be a Calabi--Yau 3-fold over $\C$, and $\O_X(1)$ a
very ample line bundle on $X$. Our definition of Calabi--Yau 3-fold
requires $X$ to be projective, with $H^1(\O_X)=0$. Write $\coh(X)$
for the abelian category of coherent sheaves on $X$, and $K(X)$ for
the numerical Grothendieck group of $\coh(X)$. Let $\tau$ denote
Gieseker stability of coherent sheaves w.r.t.\ $\O_X(1)$. If $E$ is
a coherent sheaf on $X$ then $[E]\in K(X)$ is in effect the Chern
character $\ch(E)$ of $E$ in~$H^{\rm even}(X;\Q)$.

For $\al\in K(X)$ we can form the coarse moduli schemes
$\M_\ss^\al(\tau)$, $\M_\st^\al(\tau)$ of $\tau$-(semi)stable
sheaves $E$ with $[E]=\al$. Then $\M_\ss^\al(\tau)$ is a projective
$\C$-scheme whose points correspond to S-equivalence classes of
$\tau$-semistable sheaves, and $\M_\st^\al(\tau)$ is an open
subscheme of $\M_\ss^\al(\tau)$ whose points correspond to
isomorphism classes of $\tau$-stable sheaves.

For Chern characters $\al$ with $\M_\ss^\al(\tau)=\M_\st^\al(\tau)$,
following Donaldson and Thomas \cite[\S 3]{DoTh}, Thomas \cite{Thom}
constructed a symmetric obstruction theory on $\M_\st^\al(\tau)$ and
defined the {\it Donaldson--Thomas invariant\/} to be the virtual
class
\begin{equation*}
DT^\al(\tau)=\ts\int_{[\M_\st^\al(\tau)]^\vir}1\in\Z,
\end{equation*}
an integer which `counts' $\tau$-semistable sheaves in class $\al$.
Thomas' main result \cite[\S 3]{Thom} is that $DT^\al(\tau)$ is
unchanged under deformations of the underlying Calabi--Yau 3-fold
$X$. Later, Behrend \cite{Behr} showed that Donaldson--Thomas
invariants can be written as a weighted Euler characteristic
\begin{equation*}
DT^\al(\tau)=\chi\bigl(\M_\st^\al(\tau),\nu_{\M_\st^\al(\tau)}\bigr),
\end{equation*}
where $\nu_{\M_\st^\al(\tau)}$ is the {\it Behrend function}, a
constructible function on $\M_\st^\al(\tau)$ depending only on
$\M_\st^\al(\tau)$ as a $\C$-scheme.

Conventional Donaldson--Thomas invariants $DT^\al(\tau)$ are only
defined for classes $\al$ with $\M_\ss^\al(\tau)=\M_\st^\al(\tau)$,
that is, when there are no strictly $\tau$-semistable sheaves. Also,
although $DT^\al(\tau)$ depends on the stability condition $\tau$,
that is, on the choice of very ample line bundle $\O_X(1)$ on $X$,
this dependence was not understood until now. The main goal of
\cite{JoSo} is to address these two issues.

For a Calabi--Yau 3-fold $X$ over $\C$ we will define {\it
generalized Donaldson--Thomas invariants\/}
$\bar{DT}{}^\al(\tau)\in\Q$ for all $\al\in K(X)$, which `count'
$\tau$-semistable sheaves in class $\al$. These have the following
important properties:
\begin{itemize}
\setlength{\itemsep}{0pt}
\setlength{\parsep}{0pt}
\item $\bar{DT}{}^\al(\tau)\in\Q$ is unchanged by deformations of the
Calabi--Yau 3-fold~$X$.
\item If $\M_\ss^\al(\tau)=\M_\st^\al(\tau)$ then $\bar{DT}{}^\al(\tau)$
lies in $\Z$ and equals the conventional Donaldson--Thomas invariant
$DT^\al(\tau)$ defined by Thomas~\cite{Thom}.
\item If $\M_\ss^\al(\tau)\ne\M_\st^\al(\tau)$ then conventional
Donaldson--Thomas invariants $DT^\al(\tau)$ are not defined for
class $\al$. Our generalized invariant $\bar{DT}{}^\al(\tau)$ may
lie in $\Q$ because strictly semistable sheaves $E$ make
(complicated) $\Q$-valued contributions to $\bar{DT}{}^\al(\tau)$.
For `generic' $\tau$ we have a conjecture that writes the
$\bar{DT}{}^\al(\tau)$ in terms of other, integer-valued
invariants~$\hat{DT}{}^\al(\tau)$.
\item If $\tau,\ti\tau$ are two stability conditions on $\coh(X)$,
there is an explicit change of stability condition formula giving
$\bar{DT}{}^\al(\ti\tau)$ in terms of the $\bar{DT}{}^\be(\tau)$.
\end{itemize}
These invariants are a continuation of the author's
programme~\cite{Joyc1,Joyc2,Joyc3,Joyc4,Joyc5,Joyc6,Joyc7}.

We begin in \S\ref{ig2} with some background material on {\it
constructible functions\/} and {\it stack functions\/} on Artin
stacks, taken from \cite{Joyc1,Joyc2}. Then \S\ref{ig3} summarizes
ideas from \cite{Joyc3,Joyc4,Joyc5,Joyc6} on Euler-characteristic
type invariants $J^\al(\tau)$ counting sheaves on Calabi--Yau
3-folds and their wall-crossing under change of stability condition,
and facts on Donaldson--Thomas invariants from Thomas \cite{Thom}
and Behrend~\cite{Behr}.

Section \ref{ig4} summarizes \cite[\S 5--\S 6]{JoSo}, and is the
heart of the paper. Let $X$ be a Calabi--Yau 3-fold, and $\fM$ the
moduli stack of coherent sheaves on $X$. Write $\bar\chi:K(X)\t
K(X)\ra\Z$ for the Euler form of $\coh(X)$. We will explain that the
Behrend function $\nu_\fM$ of $\fM$ satisfies two important
identities
\begin{gather*}
\nu_{\fM}(E_1\op E_2)=(-1)^{\bar\chi([E_1],[E_2])}
\nu_{\fM}(E_1)\nu_{\fM}(E_2),\\
\begin{split}
\int_{\begin{subarray}{l}[\la]\in\mathbb{P}(\Ext^1(E_2,E_1)):\\
\la\; \Leftrightarrow\; 0\ra E_1\ra F\ra E_2\ra
0\end{subarray}}\!\!\!\!\!\!\nu_{\fM}(F) \d\chi -
\int_{\begin{subarray}{l}[\la']\in\mathbb{P}(\Ext^1(E_1,E_2)):\\
\la'\; \Leftrightarrow\; 0\ra E_2\ra F'\ra E_1\ra
0\end{subarray}}\!\!\!\!\!\!\nu_{\fM}(F')\d\chi \\
=\bigl(\dim\Ext^1(E_2,E_1)-\dim\Ext^1(E_1,E_2)\bigr)
\nu_{\fM}(E_1\op E_2).
\end{split}
\end{gather*}
We use these to define a Lie algebra morphism
$\ti\Psi:\SFai(\fM)\ra\ti L(X)$, where $\SFai(\fM)$ is a special Lie
subalgebra of the {\it Ringel--Hall algebra\/} $\SFa(\fM)$ of $X$, a
large algebra with a universal construction, and $\ti L(X)$ is a
much smaller explicit Lie algebra, the $\Q$-vector space with basis
$\ti\la^\al$ for $\al\in K(X)$, and Lie bracket
\begin{equation*}
[\ti\la^\al,\ti\la^\be]=(-1)^{\bar\chi(\al,\be)}
\bar\chi(\al,\be)\ti \la^{\al+\be}.
\end{equation*}

If $\tau$ is Gieseker stability in $\coh(X)$ and $\al\in K(X)$, we
define an element $\bep^\al(\tau)$ in $\SFai(\fM)$ which `counts'
$\tau$-semistable sheaves in class $\al$ in a special way. We define
the {\it generalized Donaldson--Thomas invariant\/}
$\bar{DT}{}^\al(\tau)\in\Q$ by
\begin{equation*}
\ti\Psi\bigl(\bar\ep^\al(\tau)\bigr)=-\bar{DT}{}^\al(\tau)\ti
\la^\al.
\end{equation*}
By results in \cite{Joyc6}, the $\bep^\al(\tau)$ transform according
to a universal transformation law in the Lie algebra $\SFai(\fM)$
under change of stability condition. Applying $\ti\Psi$ shows that
$-\bar{DT}{}^\al(\tau)\ti\la^\al$ transform according to the same
law in $\ti L(X)$. This yields a wall-crossing formula for two
stability conditions $\tau,\ti\tau$ on~$\coh(X)$:
\e
\begin{aligned}
&\bar{DT}{}^\al(\ti\tau)=\\
&\!\!\!\sum_{\substack{\text{iso.}\\ \text{classes}\\
\text{of finite}\\ \text{sets $I$}}}\,
\sum_{\substack{\ka:I\ra C(X):\\ \sum\limits_{i\in I}\ka(i)=\al}}\,\,
\sum_{\begin{subarray}{l} \text{connected,}\\
\text{simply-}\\ \text{connected}\\ \text{digraphs $\Ga$,}\\
\text{vertices $I$}\end{subarray}}
\begin{aligned}[t]
(-1)^{\md{I}-1} V(I,\Ga,\ka;\tau,\ti\tau)\cdot
\prod\nolimits_{i\in I} \bar{DT}{}^{\ka(i)}(\tau)&\\
\cdot (-1)^{\frac{1}{2}\sum_{i,j\in
I}\md{\bar\chi(\ka(i),\ka(j))}}\cdot\! \prod\limits_{\text{edges
\smash{$\mathop{\bu}\limits^{\sst i}\ra\mathop{\bu}\limits^{\sst
j}$} in $\Ga$}\!\!\!\!\!\!\!\!\!\!\!\!\!\!\!\!\!\!\!\!\!\!\!
\!\!\!\!\!\!\!\!} \bar\chi(\ka(i),\ka(j))&,
\end{aligned}
\end{aligned}
\label{ig1eq}
\e
where $V(I,\Ga,\ka;\tau,\ti\tau)\in\Q$ are combinatorial
coefficients, and there are only finitely many nonzero terms.

To prove that $\bar{DT}{}^\al(\tau)$ is unchanged under deformations
of $X$, we introduce auxiliary invariants $PI^{\al,n}(\tau')\in\Z$
counting `stable pairs' $s:\O(-n)\ra E$, for $n\gg 0$ and
$E\in\coh(X)$ $\tau$-semistable in class $\al\in K(X)$. The moduli
space $\M_\stp^{\al,n}(\tau')$ of such stable pairs is a proper fine
moduli $\C$-scheme with a symmetric obstruction theory, so by the
same proof as for Donaldson--Thomas invariants \cite{Thom}, the
virtual count $PI^{\al,n}(\tau')$ of $\M_\stp^{\al,n}(\tau')$ is
deformation-invariant. By a wall-crossing proof similar to that for
\eq{ig1eq} we find that
\begin{equation*}
PI^{\al,n}(\tau')=\!\!\!\!\!\!\!\!\!\!\!\!
\sum_{\begin{subarray}{l} \al_1,\ldots,\al_l\in
C(X),\\ l\ge 1:\; \al_1
+\cdots+\al_l=\al,\\
\tau(\al_i)=\tau(\al),\text{ all\/ $i$}
\end{subarray} \!\!\!\!\!\!\!\!\! }
\begin{aligned}[t] \frac{(-1)^l}{l!} &\prod_{i=1}^{l}\bigl[
(-1)^{\bar\chi([\O_X(-n)]-\al_1-\cdots-\al_{i-1},\al_i)} \\
&\bar\chi\bigl([\O_{X}(-n)]\!-\!\al_1\!-\!\cdots\!-\!\al_{i-1},\al_i
\bigr) \bar{DT}{}^{\al_i}(\tau)\bigr].\!\!\!\!\!\!\!\!\!\!
\end{aligned}
\end{equation*}
Using deformation-invariance of the $PI^{\al,n}(\tau')$ and
induction on $\rank\al$ we find that $\bar{DT}{}^\al(\tau)$ is
deformation-invariant.

Examples show that in general the $\bar{DT}{}^\al(\tau)$ lie in $\Q$
rather than $\Z$. So it is an interesting question whether we can
rewrite the $\bar{DT}{}^\al(\tau)$ in terms of some system of
$\Z$-valued invariants, just as $\Q$-valued Gromov--Witten
invariants of Calabi--Yau 3-folds are (conjecturally) written in
terms of $\Z$-valued Gopakumar--Vafa invariants \cite{GoVa}. We
define new {\it BPS invariants\/} $\hat{DT}{}^\al(\tau)$ for $\al\in
C(X)$ to satisfy
\begin{equation*}
\bar{DT}{}^\al(\tau)=\sum_{m\ge 1,\; m\mid\al}\frac{1}{m^2}\,
\hat{DT}{}^{\al/m}(\tau),
\end{equation*}
and we conjecture that $\hat{DT}{}^\al(\tau)\in\Z$ for all $\al$ if
the stability condition $\tau$ is `generic'. Evidence for this
conjecture is given in~\cite[\S 6.1--\S 6.5 \& \S 7.6]{JoSo}.

Section \ref{ig5} summarizes \cite[\S 7]{JoSo}, which develops an
analogue of Donaldson--Thomas theory for representations of quivers
with relations coming from a superpotential. This provides a kind of
toy model for Donaldson--Thomas invariants using only polynomials
and finite-dimensional algebra, and is a source of many simple,
explicit examples. Counting invariants for quivers with
superpotential have been studied by Nakajima, Reineke, Szendr\H oi
and other authors for some years \cite{EnRe,Naga,NaNa,Naka,
MoRe,Rein1,Rein2,Szen}, under the general name of `noncommutative
Donaldson--Thomas invariants'. Curiously, the invariants studied so
far are the analogues of our pair invariants $PI^{\al,n}(\tau')$,
and the analogues of $\bar{DT}{}^\al(\tau),\hat{DT}{}^\al(\tau)$
seem to have received no attention, although they appear to the
author to be more fundamental.

A recent paper by Kontsevich and Soibelman \cite{KoSo1}, summarized
in \cite{KoSo2}, has considerable overlap with both \cite{JoSo} and
the already published
\cite{Joyc1,Joyc2,Joyc3,Joyc4,Joyc5,Joyc6,Joyc7}. The two were
completed largely independently, and the first versions of
\cite{JoSo,KoSo1} appeared on the arXiv within a few days of each
other. Kontsevich and Soibelman are far more ambitious than us,
working in triangulated categories rather than abelian categories,
over general fields $\K$ rather than $\C$, and with general motivic
invariants rather than the Euler characteristic. But for this
reason, almost every major result in \cite{KoSo1} depends explicitly
or implicitly on conjectures. The author would like to acknowledge
the contribution of \cite{KoSo1} to the ideas on
$\hat{DT}{}^\al(\tau)$ and integrality in \S\ref{ig44} below, and to
the material on quivers with superpotential in \S\ref{ig5}. The
relationship between \cite{JoSo} and \cite{KoSo1} is discussed in
detail in~\cite[\S 1.6]{JoSo}.
\medskip

\noindent{\it Acknowledgements.} The author would like to thank Tom
Bridgeland, Richard Thomas, Bal\'{a}zs Szendr\H{o}i, and his
co-author Yinan Song. This research was supported by EPSRC
grant~EP/D077990/1.

\section{Constructible functions and stack functions}
\label{ig2}

We begin with some background material on Artin stacks,
constructible functions, and stack functions, drawn from
\cite{Joyc1,Joyc2}. We restrict to the field~$\K=\C$.

\subsection{Artin stacks and constructible functions}
\label{ig21}

{\it Artin stacks\/} are a class of geometric spaces, generalizing
schemes and algebraic spaces. For a good introduction to Artin
stacks see G\'omez \cite{Gome}, and for a thorough treatment see
Laumon and Moret-Bailly \cite{LaMo}. We work throughout over the
field $\C$. We make the convention that all Artin stacks in this
paper are {\it locally of finite type}, with {\it affine geometric
stabilizers}, that is, all stabilizer groups $\Iso_\fF(x)$ are
affine algebraic $\C$-groups, and substacks are {\it locally
closed}.

Artin $\C$-stacks form a 2-{\it category}. That is, we have {\it
objects} which are $\C$-stacks $\fF,\fG$, and also two kinds of
morphisms, 1-{\it morphisms} $\phi,\psi:\fF\ra\fG$ between
$\C$-stacks, and 2-{\it morphisms} $A:\phi\ra\psi$ between
1-morphisms.

\begin{dfn} Let $\fF$ be a $\C$-stack. Write $\fF(\C)$ for the set
of 2-isomorphism classes $[x]$ of 1-morphisms $x:\Spec\C\ra\fF$.
Elements of $\fF(\C)$ are called $\C$-{\it points\/} of $\fF$. If
$\phi:\fF\ra\fG$ is a 1-morphism then composition with $\phi$
induces a map of sets~$\phi_*:\fF(\C)\ra\fG(\C)$.

For a 1-morphism $x:\Spec\C\ra\fF$, the {\it stabilizer group}
$\Iso_\fF(x)$ is the group of 2-morphisms $A:x\ra x$. When $\fF$ is
an Artin $\C$-stack, $\Iso_\fF(x)$ is an {\it algebraic $\C$-group},
which we assume is affine. If $\phi:\fF\ra\fG$ is a 1-morphism,
composition induces a morphism of $\C$-groups
$\phi_*:\Iso_\fF([x])\ra\Iso_\fG\bigr(\phi_*([x])\bigr)$,
for~$[x]\in\fF(\C)$.
\label{ig2def1}
\end{dfn}

We discuss {\it constructible functions} on $\C$-stacks,
following~\cite{Joyc1}.

\begin{dfn} Let $\fF$ be an Artin $\C$-stack. We call
$C\subseteq\fF(\C)$ {\it constructible} if $C=\bigcup_{i\in I}
\fF_i(\C)$, where $\{\fF_i:i\in I\}$ is a finite collection of
finite type Artin $\C$-substacks $\fF_i$ of $\fF$. We call
$S\subseteq\fF(\C)$ {\it locally constructible} if $S\cap C$ is
constructible for all constructible $C\subseteq\fF(\C)$. A function
$f:\fF(\C)\ra\Q$ is called {\it constructible} if $f(\fF(\C))$ is
finite and $f^{-1}(c)$ is a constructible set in $\fF(\C)$ for each
$c\in f(\fF(\C))\sm\{0\}$. A function $f:\fF(\C)\ra\Q$ is called
{\it locally constructible} if $f\cdot\de_C$ is constructible for
all constructible $C\subseteq\fF(\C)$, where $\de_C$ is the
characteristic function of $C$. Write $\CF(\fF)$ and $\LCF(\fF)$ for
the $\Q$-vector spaces of $\Q$-valued constructible and locally
constructible functions on~$\fF$.
\label{ig2def2}
\end{dfn}

Following \cite[\S 4--\S 5]{Joyc1} we define {\it pushforwards} and
{\it pullbacks} of constructible functions along 1-morphisms.

\begin{dfn} Let $\fF,\fG$ be Artin $\C$-stacks and $\phi:\fF\ra\fG$
a representable 1-morphism. For $f\in\CF(\fF)$, define
$\CF^\stk(\phi)f:\fG(\C)\ra\Q$ by
\begin{equation*}
\CF^\stk(\phi)f(y)=\chi\bigl(\fF\t_{\phi,\fG,y}\Spec\C,\pi_\fF^*(f)\bigr)
\quad\text{for $y\in\fG(\C)$,}
\end{equation*}
where $\fF\t_{\phi,\fG,y}\Spec\C$ is a $\C$-scheme (or algebraic
space) as $\phi$ is representable, and $\chi(\cdots)$ is the Euler
characteristic of this $\C$-scheme weighted by $\pi_\fF^*(f)$. Then
$\CF^\stk(\phi):\CF(\fF)\ra\CF(\fG)$ is a $\Q$-linear map called the
{\it stack pushforward}.

Let $\th:\fF\ra\fG$ be a finite type 1-morphism. The {\it
pullback\/} $\th^*:\CF(\fG)\ra\CF(\fF)$ is given by $\th^*(f)=
f\ci\th_*$. It is a $\Q$-linear map.
\label{ig2def3}
\end{dfn}

Here \cite[\S 4--\S 5]{Joyc1} are some properties of these.

\begin{thm} Let\/ $\fE,\fF,\fG,\fH$ be Artin $\C$-stacks and\/
$\be:\fF\ra\fG,$ $\ga:\fG\ra\fH$ be $1$-morphisms. Then
\ea
\CF^\stk(\ga\ci\be)&=\CF^\stk(\ga)\ci\CF^\stk(\be):\CF(\fF)\ra\CF(\fH),
\label{ig2eq1}\\
(\ga\ci\be)^*&=\be^*\ci\ga^*:\CF(\fH)\ra\CF(\fF),
\label{ig2eq2}
\ea
supposing $\be,\ga$ representable in {\rm\eq{ig2eq1},} and of finite
type in \eq{ig2eq2}. If
\begin{equation*}
\begin{gathered}
\xymatrix{
\fE \ar[r]_\eta \ar[d]^\th & \fG \ar[d]_\psi \\
\fF \ar[r]^\phi & \fH }
\end{gathered}
\quad
\begin{gathered}
\text{is a Cartesian square with}\\
\text{$\eta,\phi$ representable and}\\
\text{$\th,\psi$ of finite type, then}\\
\text{the following commutes:}
\end{gathered}
\quad
\begin{gathered}
\xymatrix@C=35pt{
\CF(\fE) \ar[r]_{\CF^\stk(\eta)} & \CF(\fG) \\
\CF(\fF) \ar[r]^{\CF^\stk(\phi)} \ar[u]_{\th^*} & \CF(\fH).
\ar[u]^{\psi^*} }
\end{gathered}
\end{equation*}
\label{ig2thm1}
\end{thm}

\subsection{Stack functions}
\label{ig22}

{\it Stack functions\/} are a universal generalization of
constructible functions introduced in \cite[\S 3]{Joyc2}. Here
\cite[Def.~3.1]{Joyc2} is the basic definition.

\begin{dfn} Let $\fF$ be an Artin $\C$-stack. Consider pairs
$(\fR,\rho)$, where $\fR$ is a finite type Artin $\C$-stack and
$\rho:\fR\ra\fF$ is a representable 1-morphism. We call two pairs
$(\fR,\rho)$, $(\fR',\rho')$ {\it equivalent\/} if there exists a
1-isomorphism $\io:\fR\ra\fR'$ such that $\rho'\ci\io$ and $\rho$
are 2-isomorphic 1-morphisms $\fR\ra\fF$. Write $[(\fR,\rho)]$ for
the equivalence class of $(\fR,\rho)$. If $(\fR,\rho)$ is such a
pair and $\fS$ is a closed $\C$-substack of $\fR$ then
$(\fS,\rho\vert_\fS)$, $(\fR\sm\fS,\rho\vert_{\fR\sm\fS})$ are pairs
of the same kind.

Define $\SF(\fF)$ to be the $\Q$-vector space generated by
equivalence classes $[(\fR,\rho)]$ as above, with for each closed
$\C$-substack $\fS$ of $\fR$ a relation
\e
[(\fR,\rho)]=[(\fS,\rho\vert_\fS)]+[(\fR\sm\fS,\rho\vert_{\fR\sm\fS})].
\label{ig2eq3}
\e
\label{ig2def4}
\end{dfn}\vspace{-17pt}

Elements of $\SF(\fF)$ will be called {\it stack functions}. We
relate $\CF(\fF)$ and~$\SF(\fF)$.

\begin{dfn} Let $\fF$ be an Artin $\C$-stack and
$C\subseteq\fF(\C)$ be constructible. Then
$C=\coprod_{i=1}^n\fR_i(\C)$, for $\fR_1,\ldots,\fR_n$ finite type
$\C$-substacks of $\fF$. Let $\rho_i:\fR_i\ra\fF$ be the inclusion
1-morphism. Then $[(\fR_i,\rho_i)]\in\SF(\fF)$. Define
$\bde_C=\ts\sum_{i=1}^n[(\fR_i,\rho_i)]\in\SF(\fF)$. We think of
this as the analogue of the characteristic function
$\de_C\in\CF(\fF)$ of $C$. Define a $\Q$-linear map
$\io_\fF:\CF(\fF)\ra\SF(\fF)$ by $\io_\fF(f)=\ts\sum_{0\ne c\in
f(\fF(\C))}c\cdot\bde_{f^{-1}(c)}$. Define $\Q$-linear
$\pi_\fF^\stk:\SF(\fF)\ra\CF(\fF)$ by
\begin{equation*}
\pi_\fF^\stk\bigl(\ts\sum_{i=1}^nc_i[(\fR_i,\rho_i)]\bigr)=
\ts\sum_{i=1}^nc_i\CF^\stk(\rho_i)1_{\fR_i},
\end{equation*}
where $1_{\fR_i}$ is the function $1\in\CF(\fR_i)$. Then
$\pi_\fF^\stk\ci\io_\fF$ is the identity on~$\CF(\fF)$.
\label{ig2def5}
\end{dfn}

The operations on constructible functions in \S\ref{ig21} extend to
stack functions.

\begin{dfn} Let $\phi:\fF\!\ra\!\fG$ be a representable 1-morphism
of Artin $\C$-stacks. Define the {\it pushforward\/}
$\phi_*:\SF(\fF)\!\ra\!\SF(\fG)$~by
\e
\phi_*:\ts\sum_{i=1}^mc_i[(\fR_i,\rho_i)]\longmapsto
\ts\sum_{i=1}^mc_i[(\fR_i,\phi\ci\rho_i)].
\label{ig2eq4}
\e
Let $\phi:\fF\!\ra\!\fG$ be of finite type. Define the {\it
pullback\/} $\phi^*:\SF(\fG)\!\ra\!\SF(\fF)$~by
\e
\phi^*:\ts\sum_{i=1}^mc_i[(\fR_i,\rho_i)]\longmapsto
\ts\sum_{i=1}^mc_i[(\fR_i\t_{\rho_i,\fG,\phi}\fF,\pi_\fF)].
\label{ig2eq5}
\e
The {\it tensor product\/} $\ot:\SF(\fF)\t\SF(\fG)\ra
\SF(\fF\t\fG)$~is
\e
\bigl(\ts\sum_{i=1}^mc_i[(\fR_i,\rho_i)]\bigr)\!\ot\!
\bigl(\ts\sum_{j=1}^nd_j[(\fS_j,\si_j)]\bigr)\!=\!\ts
\sum_{i,j}c_id_j[(\fR_i\!\t\!\fS_j,\rho_i\!\t\!\si_j)].
\label{ig2eq6}
\e
\label{ig2def6}
\end{dfn}\vspace{-17pt}

Here \cite[Th.~3.5]{Joyc2} is the analogue of Theorem~\ref{ig2thm1}.

\begin{thm} Let\/ $\fE,\fF,\fG,\fH$ be Artin $\C$-stacks and\/
$\be:\fF\ra\fG,$ $\ga:\fG\ra\fH$ be $1$-morphisms. Then
\begin{equation*}
(\ga\!\ci\!\be)_*\!=\!\ga_*\!\ci\!\be_*:\SF(\fF)\!\ra\!\SF(\fH),\qquad
(\ga\!\ci\!\be)^*\!\!=\!\be^*\!\ci\!\ga^*:\SF(\fH)\!\ra\!\SF(\fF),
\end{equation*}
for $\be,\ga$ representable in the first equation, and of finite
type in the second. If
\begin{equation*}
\begin{gathered}
\xymatrix@R=15pt{
\fE \ar[r]_\eta \ar[d]^{\,\th} & \fG \ar[d]_{\psi\,} \\
\fF \ar[r]^\phi & \fH }
\end{gathered}
\quad
\begin{gathered}
\text{is a Cartesian square with}\\
\text{$\th,\psi$ of finite type and}\\
\text{$\eta,\phi$ representable, then}\\
\text{the following commutes:}
\end{gathered}
\quad
\begin{gathered}
\xymatrix@C=35pt@R=10pt{
\SF(\fE) \ar[r]_{\eta_*} & \SF(\fG) \\
\SF(\fF) \ar[r]^{\phi_*} \ar[u]_{\,\th^*} & \SF(\fH).
\ar[u]^{\psi^*\,} }
\end{gathered}
\end{equation*}
\label{ig2thm2}
\end{thm}

In \cite[\S 3]{Joyc2} we relate pushforwards and pullbacks of stack
and constructible functions using~$\io_\fF,\pi_\fF^\stk$.

\begin{thm} Let\/ $\phi:\fF\ra\fG$ be a $1$-morphism of Artin
$\C$-stacks. Then
\begin{itemize}
\setlength{\itemsep}{0pt}
\setlength{\parsep}{0pt}
\item[{\rm(a)}] $\phi^*\!\ci\!\io_\fG\!=\!\io_\fF\!\ci\!\phi^*:
\CF(\fG)\!\ra\!\SF(\fF)$ if\/ $\phi$ is of finite type;
\item[{\rm(b)}] $\pi^\stk_\fG\ci\phi_*=\CF^\stk(\phi)\ci\pi_\fF^\stk:
\SF(\fF)\ra\CF(\fG)$ if\/ $\phi$ is representable; and
\item[{\rm(c)}] $\pi^\stk_\fF\ci\phi^*=\phi^*\ci\pi_\fG^\stk:
\SF(\fG)\ra\CF(\fF)$ if\/ $\phi$ is of finite type.
\end{itemize}
\label{ig2thm3}
\end{thm}

We define some projections $\Pi^\vi_n:\SF(\fF)\ra\SF(\fF)$, \cite[\S
5]{Joyc2}.

\begin{dfn} For any Artin $\C$-stack $\fF$ we will define linear maps
$\Pi^\vi_n:\SF(\fF)\ra\SF(\fF)$ for $n\ge 0$. Now $\SF(\fF)$ is
generated by $[(\fR,\rho)]$ with $\fR$ 1-isomorphic to a quotient
$[X/G]$, for $X$ a quasiprojective $\C$-variety and $G$ a special
algebraic $\C$-group, with maximal torus~$T^G$.

Let ${\cal S}(T^G)$ be the set of subsets of $T^G$ defined by
Boolean operations upon closed $\C$-subgroups $L$ of $T^G$. Define a
measure $\d\mu_n:{\cal S}(T^G)\ra\Z$ to be additive upon disjoint
unions of sets in ${\cal S}(T^G)$, and to satisfy $\d\mu_n(L)=1$ if
$\dim L=n$ and $\d\mu_n(L)=0$ if $\dim L\ne 0$ for all algebraic
$\C$-subgroups $L$ of $T^G$. Define
\e
\begin{split}
&\Pi^\vi_n\bigl([(\fR,\rho)]\bigr)=\\
&\int_{t\in T^G}\frac{\md{\{w\in W(G,T^G):w\cdot
t=t\}}}{\md{W(G,T^G)}}\,\bigl[\bigl([X^{\{t\}}/
C_G(\{t\})],\rho\ci\io^{\{t\}}\bigr)\bigr]\d\mu_n.
\end{split}
\label{ig2eq7}
\e
Here $X^{\{t\}}$ is the subscheme of $X$ fixed by $t$, and
$C_G(\{t\})$ is the centralizer of $t$ in $G$, and
$\io^{\{t\}}:[X^{\{t\}}/ C_G(\{t\})]\ra[X/G]$ is the obvious
1-morphism.

The integrand in \eq{ig2eq7}, regarded as a function of $t\in T^G$,
is a constructible function taking only finitely many values. The
level sets of the function lie in ${\cal S}(T^G)$, so they are
measurable w.r.t.\ $\d\mu_n$, and the integral is well-defined. In
\cite[\S 5]{Joyc2} we show \eq{ig2eq7} induces a unique linear
map~$\Pi^\vi_n:\SF(\fF)\ra\SF(\fF)$.
\label{ig2def7}
\end{dfn}

Here \cite[\S 5]{Joyc2} are some properties of the $\Pi^\vi_n$.

\begin{thm} In the situation above, we have:
\begin{itemize}
\setlength{\itemsep}{0pt}
\setlength{\parsep}{0pt}
\item[{\rm(i)}] $(\Pi^\vi_n)^2=\Pi^\vi_n,$ so that\/ $\Pi^\vi_n$ is
a projection, and\/ $\Pi^\vi_m\ci\Pi^\vi_n=0$ for~$m\ne n$.
\item[{\rm(ii)}] For all\/ $f\in\SF(\fF)$ we have $f=\sum_{n\ge 0}
\Pi^\vi_n(f),$ where the sum makes sense as $\Pi^\vi_n(f)=0$
for~$n\gg 0$.
\item[{\rm(iii)}] If\/ $\phi:\fF\ra\fG$ is a $1$-morphism of Artin
$\C$-stacks then~$\Pi^\vi_n\ci\phi_*=\phi_*\ci\Pi^\vi_n:\SF(\fF)\ra
\SF(\fG)$.
\item[{\rm(iv)}] If\/ $f\in\SF(\fF),$ $g\in\SF(\fG)$
then~$\Pi^\vi_n(f\ot
g)=\sum_{m=0}^n\Pi^\vi_m(f)\ot\Pi^\vi_{n-m}(g)$.
\end{itemize}
\label{ig2thm4}
\end{thm}

Roughly speaking, $\Pi^\vi_n$ projects $[(\fR,\rho)]\in\SF(\fF)$ to
$[(\fR_n,\rho)]$, where $\fR_n$ is the substack of points
$r\in\fR(\C)$ whose stabilizer groups $\Iso_\fR(r)$ have rank~$n$.

\subsection{Stack function spaces
$\bar{\rm SF}({\mathfrak F},\chi,{\mathbb Q})$}
\label{ig23}

We will also need another family of spaces $\oSF(\fF,\chi,\Q)$, from
\cite[\S 5--\S 6]{Joyc2}.

\begin{dfn} Let $\fF$ be an Artin $\C$-stack. Consider pairs
$(\fR,\rho)$, where $\fR$ is a finite type Artin $\C$-stack and
$\rho:\fR\ra\fF$ is a representable 1-morphism, with equivalence as
in Definition \ref{ig2def4}. Define $\oSF(\fF,\chi,\Q)$ to be the
$\Q$-vector space generated by equivalence classes $[(\fR,\rho)]$,
with the following relations:
\begin{itemize}
\setlength{\itemsep}{0pt}
\setlength{\parsep}{0pt}
\item[(i)] Given $[(\fR,\rho)]$ as above and $\fS$ a closed $\C$-substack
of $\fR$ we have $[(\fR,\rho)]=[(\fS,\rho\vert_\fS)]+[(\fR\sm\fS,
\rho\vert_{\fR\sm\fS})]$, as in~\eq{ig2eq3}.
\item[(ii)] Let $\fR$ be a finite type Artin $\C$-stack, $U$ a
quasiprojective $\C$-variety, $\pi_\fR:\fR\t U\ra\fR$ the natural
projection, and $\rho:\fR\ra\fF$ a 1-morphism. Then~$[(\fR\t
U,\rho\ci\pi_\fR)] =\chi([U])[(\fR,\rho)]$.

Here $\chi(U)\in\Z$ is the {\it Euler characteristic\/} of $U$. It
is a {\it motivic invariant\/} of $\C$-schemes, that is,
$\chi(U)=\chi(V)+\chi(U\sm V)$ for $V\subset U$ closed.
\item[(iii)] Given $[(\fR,\rho)]$ as above and a 1-isomorphism
$\fR\cong[X/G]$ for $X$ a quasiprojective $\C$-variety and $G$ a
very special algebraic $\C$-group acting on $X$ with maximal torus
$T^G$, we have
\begin{equation*}
[(\fR,\rho)]=\ts\sum_{Q\in\cQ(G,T^G)}F(G,T^G,Q)
\bigl[\bigl([X/Q],\rho\ci\io^Q\bigr)\bigr],
\end{equation*}
where $\io^Q:[X/Q]\ra\fR\cong[X/G]$ is the natural projection
1-morphism.
\end{itemize}
Here $\cQ(G,T^G)$ is a certain finite set of $\C$-subgroups of
$T^G$, and $F(G,T^G,Q)\in\Q$ are a system of rational coefficients
defined in \cite[\S 6.2]{Joyc2}. Define
$\bar\Pi^{\chi,\Q}_\fF:\SF(\fF) \ra\oSF(\fF,\chi,\Q)$ by
$\bar\Pi^{\chi,\Q}_\fF:\ts\sum_{i\in I}c_i[(\fR_i,\rho_i)]\mapsto
\ts\sum_{i\in I}c_i[(\fR_i,\rho_i)]$. Define  {\it pushforwards\/}
$\phi_*$, {\it pullbacks\/} $\phi^*$, {\it tensor products} $\ot$
and {\it projections\/} $\Pi^\vi_n$ on the spaces $\oSF(*,\chi,\Q)$
as in \S\ref{ig22}. The important point is that
\eq{ig2eq4}--\eq{ig2eq7} are compatible with the relations defining
$\oSF(*,\chi,\Q)$, or they would not be well-defined. The analogues
of Theorems \ref{ig2thm2}, \ref{ig2thm3} and \ref{ig2thm4} hold
for~$\oSF(*,\chi,\Q)$.
\label{ig2def8}
\end{dfn}

Here \cite[\S 5--\S 6]{Joyc2} is a useful way to represent these
spaces. It means that by working in $\oSF(\fF,\chi,\Q),$ we can
treat all stabilizer groups as if they are abelian.

\begin{prop} $\oSF(\fF,\chi,\Q)$ is spanned over $\Q$ by
$[(U\t[\Spec\C/T],\rho)],$ for\/ $U$ a quasiprojective $\C$-variety
and\/ $T$ an algebraic $\C$-group isomorphic to $\bG_m^k\t K$ for
$k\ge 0$ and\/ $K$ finite abelian. Moreover
\begin{equation*}
\Pi^\vi_n\bigl([(U\t[\Spec\C/T],\rho)]\bigr)=
\begin{cases} [(U\t[\Spec\C/T],\rho)], &\dim T=n,\\
0, & \text{otherwise.}\end{cases}
\end{equation*}
\label{ig2prop}
\end{prop}

\section{Background material on Calabi--Yau 3-folds}
\label{ig3}

We now summarize some facts on Donaldson--Thomas invariants and
other sheaf-counting invariants on Calabi--Yau 3-folds prior to our
book \cite{JoSo}. Sections \ref{ig31}--\ref{ig33} review material
from the author's series of papers \cite{Joyc3,Joyc4,Joyc5,Joyc6},
and \S\ref{ig34} explains results on Donaldson--Thomas theory from
Thomas \cite{Thom} and Behrend \cite{Behr}. For simplicity we
restrict to Calabi--Yau 3-folds and to the field $\K=\C$, although
much of \cite{Behr,Joyc3,Joyc4,Joyc5,Joyc6,Thom} works in greater
generality.

\subsection{The Ringel--Hall algebra of a Calabi--Yau 3-fold}
\label{ig31}

We will use the following notation for the rest of the paper.

\begin{dfn} A {\it Calabi--Yau\/ $3$-fold\/} is a smooth projective
3-fold $X$ over $\C$, with trivial canonical bundle $K_X$. In
\S\ref{ig4} we will also assume that $H^1(\O_X)=0$. The {\it
Grothendieck group\/} $K_0(X)$ of $\coh(X)$ is the abelian group
generated by all isomorphism classes $[E]$ of objects $E$ in
$\coh(X)$, with the relations $[E]+[G]=[F]$ for each short exact
sequence $0\ra E\ra F\ra G\ra 0$. The {\it Euler form\/}
$\bar\chi:K_0(X)\t K_0(X)\ra\Z$ is a biadditive map satisfying
\e
\bar\chi\bigl([E],[F]\bigr)=\ts\sum_{i\ge 0}(-1)^i\dim\Ext^i(E,F)
\label{ig3eq1}
\e
for all $E,F\in\coh(X)$. As $X$ is a Calabi--Yau 3-fold, Serre
duality gives $\Ext^i(F,E)\cong\Ext^{3-i}(E,F)^*$, so
$\dim\Ext^i(F,E)=\dim\Ext^{3-i}(E,F)$ for all $E,F\in\coh(X)$.
Therefore $\bar\chi$ is also given by
\e
\begin{split}
\bar\chi\bigl([E],[F]\bigr)=\,&\bigl(\dim\Hom(E,F)-\dim\Ext^1(E,F)
\bigr)-\\
&\bigl(\dim\Hom(F,E)-\dim\Ext^1(F,E)\bigr).
\end{split}
\label{ig3eq2}
\e
Thus the Euler form $\bar\chi$ on $K_0(X)$ is {\it antisymmetric}.

The {\it numerical Grothendieck group\/} $K(X)$ is the quotient of
$K_0(X)$ by the kernel of $\bar\chi$. Then $\bar\chi$ on $K_0(X)$
descends to a nondegenerate, biadditive Euler form~$\bar\chi:K(X)\t
K(X)\ra\Z$.

Define the `positive cone' $C(X)$ in $K(X)$ to be
\begin{equation*}
C(X)=\bigl\{[E]\in K(X):0\not\cong E\in\coh(X)\bigr\}\subset K(X).
\end{equation*}

Write $\fM$ for the moduli stack of objects in $\coh(X)$. It is an
Artin $\C$-stack, locally of finite type. Points of $\fM(\C)$
correspond to isomorphism classes $[E]$ of objects $E$ in $\coh(X)$,
and the stabilizer group $\Iso_{\fM}([E])$ in $\fM$ is isomorphic as
an algebraic $\C$-group to the automorphism group $\Aut(E)$. For
$\al\in C(X)$, write $\fM^\al$ for the substack of objects
$E\in\coh(X)$ in class $\al$ in $K(X)$. It is an open and closed
$\C$-substack of $\fM$.

Write $\fExact$ for the moduli stack of short exact sequences $0\ra
E_1\ra E_2\ra E_3\ra 0$ in $\coh(X)$. It is an Artin $\C$-stack,
locally of finite type. For $j=1,2,3$ write $\pi_j:\fExact\ra\fM$
for the 1-morphism projecting $0\ra E_1\ra E_2\ra E_3\ra 0$ to
$E_j$. Then $\pi_2$ is {\it representable}, and
$\pi_1\t\pi_3:\fExact\ra\fM\t\fM$ is of {\it finite type}.
\label{ig3def1}
\end{dfn}

In \cite{Joyc4} we define {\it Ringel--Hall algebras}, using stack
functions.

\begin{dfn} Define bilinear operations $*$ on
$\SF(\fM),\oSF(\fM,\chi,\Q)$ by
\begin{equation*}
f*g=(\pi_2)_*\bigl((\pi_1\t\pi_3)^*(f\ot g)\bigr),
\end{equation*}
using pushforwards, pullbacks and tensor products in Definition
\ref{ig2def6}. They are well-defined as $\pi_2$ is representable,
and $\pi_1\t\pi_3$ is of finite type. By \cite[Th.~5.2]{Joyc4},
whose proof uses Theorem \ref{ig2thm2}, this * is {\it associative},
and makes $\SF(\fM)$, $\oSF(\fM,\chi,\Q)$ into noncommutative
$\Q$-algebras, called {\it Ringel--Hall algebras}, with identity
$\bde_{[0]}$, where $[0]\in\fM$ is the zero object. The projection
$\bar\Pi^{\chi,\Q}_{\fM}:\SF(\fM)\ra\oSF(\fM,\chi,\Q)$ is an algebra
morphism.

As these algebras are inconveniently large for some purposes, in
\cite[Def.~5.5]{Joyc4} we define subalgebras $\SFa(\fM),
\oSFa(\fM,\chi,\Q)$ using the algebra structure on stabilizer groups
in $\fM$. Suppose $[(\fR,\rho)]$ is a generator of $\SF(\fM)$. Let
$r\in\fR(\C)$ with $\rho_*(r)=[E]\in\fM(\C)$, for some
$E\in\coh(X)$. Then $\rho$ induces a morphism of stabilizer
$\C$-groups $\rho_*:\Iso_\fR(r)\ra\Iso_\fM([E])\cong\Aut(E)$. As
$\rho$ is representable this is {\it injective}, and induces an
isomorphism of $\Iso_\fR(r)$ with a $\C$-subgroup of $\Aut(E)$. Now
$\Aut(E)=\End(E)^\t$ is the $\C$-group of invertible elements in a
{\it finite-dimensional\/ $\C$-algebra} $\End(E)=\Hom(E,E)$. We say
that $[(\fR,\rho)]$ {\it has algebra stabilizers\/} if whenever
$r\in\fR(\C)$ with $\rho_*(r)=[E]$, the $\C$-subgroup
$\rho_*\bigl(\Iso_\fR(r)\bigr)$ in $\Aut(E)$ is the $\C$-group
$A^\t$ of invertible elements in a $\C$-subalgebra $A$ in $\End(E)$.
Write $\SFa(\fM),\oSFa(\fM,\chi,\Q)$ for the subspaces of
$\SF(\fM),\oSF(\fM,\chi,\Q)$ spanned over $\Q$ by $[(\fR,\rho)]$
with algebra stabilizers. Then \cite[Prop.~5.7]{Joyc4} shows that
$\SFa(\fM),\oSFa(\fM,\chi,\Q)$ are {\it subalgebras\/} of the
Ringel--Hall algebras~$\SF(\fM),\oSF(\fM,\chi,\Q)$.

Now \cite[Cor.~5.10]{Joyc4} shows that $\SFa(\fM),\oSFa
(\fM,\chi,\Q)$ are closed under the operators $\Pi^\vi_n$ on
$\SF(\fM),\oSF(\fM,\chi,\Q)$ defined in \S\ref{ig22}. In
\cite[Def.~5.14]{Joyc4} we define $\SFai(\fM),\oSFai(\fM,\chi,\Q)$
to be the subspaces of $f$ in $\SFa(\fM)$ and $\oSFa (\fM,\chi,\Q)$
with $\Pi^\vi_1(f)=f$. We think of $\SFai(\fM),\oSFai(\fM,\chi,\Q)$
as stack functions {\it `supported on virtual indecomposables'}.

In \cite[Th.~5.18]{Joyc4} we show that $\SFai(\fM),\oSFai
(\fM,\chi,\Q)$ are closed under the Lie bracket $[f,g]=f*g-g*f$ on
$\SFa(\fM),\oSFa(\fM,\chi,\Q)$. Thus, $\SFai\!(\fM),\oSFai
(\fM,\chi,\Q)$ are {\it Lie subalgebras} of
$\SFa(\fM),\oSFa(\fM,\chi,\Q)$.
\label{ig3def2}
\end{dfn}

As in \cite[Cor.~5.11]{Joyc4}, Proposition \ref{ig2prop} simplifies
to give:

\begin{prop} $\oSFa(\fM,\chi,\Q)$ is spanned over $\Q$ by
elements of the form $[(U\t[\Spec\C/\bG_m^k],\rho)]$ with algebra
stabilizers, for\/ $U$ a quasiprojective $\C$-variety and\/ $k\ge
0$. Also $\oSFai(\fM,\chi,\Q)$ is spanned over $\Q$ by
$[(U\t[\Spec\C/\bG_m],\rho)]$ with algebra stabilizers, for\/ $U$ a
quasiprojective $\C$-variety.
\label{ig3prop1}
\end{prop}

All the above except \eq{ig3eq2} works for $X$ an arbitrary smooth
projective $\C$-scheme, but our next result uses the Calabi--Yau
3-fold assumption on $X$ in an essential way. We follow \cite[\S
6.5--\S 6.6]{Joyc4}, but use the notation of~\cite[\S 3.4]{JoSo}.

\begin{dfn} Define an explicit Lie algebra $L(X)$ over $\Q$ to be
the $\Q$-vector space with basis of symbols $\la^\al$ for $\al\in
K(X)$, with Lie bracket
\e
[\la^\al,\la^\be]=\bar\chi(\al,\be)\la^{\al+\be}
\label{ig3eq3}
\e
for $\al,\be\in K(X)$. As $\bar\chi$ is antisymmetric, \eq{ig3eq3}
satisfies the Jacobi identity and makes $L(X)$ into an
infinite-dimensional Lie algebra over $\Q$.

Define a $\Q$-linear map $\Psi^{\chi,\Q}:\oSFai(\fM,\chi,\Q)\ra
L(X)$ by
\e
\Psi^{\chi,\Q}(f)=\ts\sum_{\al\in K(X)}\ga^\al \la^{\al},
\label{ig3eq4}
\e
where $\ga^\al\in\Q$ is defined as follows. Proposition
\ref{ig3prop1} says $\oSFai(\fM,\chi,\Q)$ is spanned by elements
$[(U\t[\Spec\C/\bG_m],\rho)]$. We may write
\e
f\vert_{\fM^\al}=\ts\sum_{i=1}^n\de_i[(U_i\t[\Spec\C/
\bG_m],\rho_i)],
\label{ig3eq5}
\e
where $\de_i\in\Q$ and $U_i$ is a quasiprojective $\C$-variety. We
set
\begin{equation*}
\ga^\al=\ts\sum_{i=1}^n\de_i\chi(U_i).
\end{equation*}
This is independent of the choices in \eq{ig3eq5}. Now define
$\Psi:\SFai(\fM)\ra L(X)$ by~$\Psi=\Psi^{\chi,\Q}\ci\bar
\Pi^{\chi,\Q}_{\fM}$.
\label{ig3def3}
\end{dfn}

In \cite[Th.~6.12]{Joyc4}, using equation \eq{ig3eq2}, we prove:

\begin{thm} $\Psi:\SFai(\fM)\ra L(X)$ and\/ $\Psi^{\chi,\Q}:
\oSFai(\fM,\chi,\Q)\ab\ra L(X)$ are Lie algebra morphisms.
\label{ig3thm1}
\end{thm}

\subsection{Stability conditions on $\coh(X)$ and invariants
$J^\al(\tau)$}
\label{ig32}

Next we discuss material in \cite{Joyc5} on {\it stability
conditions}. We continue to use the notation of \S\ref{ig31}, with
$X$ a Calabi--Yau 3-fold.

\begin{dfn} Suppose $(T,\le)$ is a totally ordered set, and
$\tau:C(X)\ra T$ a map. We call $(\tau,T,\le)$ a {\it stability
condition} on $\coh(X)$ if whenever $\al,\be,\ga\in C(X)$ with
$\be=\al+\ga$ then either $\tau(\al)\!<\!\tau(\be)\!<\!\tau(\ga)$,
or $\tau(\al)\!>\!\tau(\be)\!>\!\tau(\ga)$, or
$\tau(\al)\!=\!\tau(\be)\!=\!\tau(\ga)$. We call $(\tau,T,\le)$ a
{\it weak stability condition} on $\coh(X)$ if whenever $\al,\be,
\ga\in C(X)$ with $\be=\al+\ga$ then either $\tau(\al)\!\le\!
\tau(\be)\!\le\!\tau(\ga)$, or $\tau(\al)\!\ge\!\tau(\be)\!\ge
\!\tau(\ga)$. For such $(\tau,T,\le)$, we call a nonzero sheaf $E$
in $\coh(X)$
\begin{itemize}
\setlength{\itemsep}{0pt}
\setlength{\parsep}{0pt}
\item[(i)] $\tau$-{\it stable} if for all $S\subset E$ with
$S\not\cong 0,E$ we have $\tau([S])<\tau([E/S])$; and
\item[(ii)] $\tau$-{\it semistable} if for all $S\subset E$ with
$S\not\cong 0,E$ we have $\tau([S])\le\tau([E/S])$.
\end{itemize}

For $\al\in C(X)$, write $\fM_\ss^\al(\tau),\fM_\st^\al(\tau)$ for
the moduli stacks of $\tau$-(semi)stable $E\in\A$ with class
$[E]=\al$ in $K(X)$. They are open $\C$-substacks of $\fM^\al$. We
call $(\tau,T,\le)$ {\it permissible} if:
\begin{itemize}
\setlength{\itemsep}{0pt}
\setlength{\parsep}{0pt}
\item[(a)] $\coh(X)$ is $\tau$-{\it artinian}, that is, there exist
no infinite chains of subobjects
$\cdots\!\subsetneq\!E_2\!\subsetneq\! E_1\!\subsetneq\!E_0=X$ in
$\A$ and $\tau([E_{n+1}])\!\ge\! \tau([E_n/E_{n+1}])$ for all $n$;
and
\item[(b)] $\fM_\ss^\al(\tau)$ is a {\it finite type\/} substack of
$\fM^\al$ for all~$\al\in C(X)$.
\end{itemize}
\label{ig3def4}
\end{dfn}

Here are two important examples:

\begin{ex} Define $G$ to be the set of monic rational polynomials in
$t$ of degree at most~3:
\begin{equation*}
G=\bigl\{p(t)=t^d+a_{d-1}t^{d-1}+\cdots+a_0:d=0,1,2,3,\;\>
a_0,\ldots,a_{d-1}\in\Q\bigr\}.
\end{equation*}
Define a {\it total order\/} `$\le$' on $G$ by $p\le p'$ for
$p,p'\in G$ if either
\begin{itemize}
\setlength{\itemsep}{0pt}
\setlength{\parsep}{0pt}
\item[(a)] $\deg p>\deg p'$, or
\item[(b)] $\deg p=\deg p'$ and $p(t)\le p'(t)$ for all $t\gg 0$.
\end{itemize}
We write $p< q$ if $p\le q$ and $p\ne q$.

Fix a very ample line bundle $\O_X(1)$ on $X$. For $E\in\coh(X)$,
the {\it Hilbert polynomial\/} $P_E$ is the unique polynomial in
$\Q[t]$ such that $P_E(n)=\dim H^0(E(n))$ for all $n\gg 0$.
Equivalently, $P_E(n)=\bar\chi\bigl([\O_X(-n)],[E]\bigr)$ for all
$n\in\Z$. Thus, $P_E$ depends only on the class $\al\in K(X)$ of
$E$, and we may write $P_\al$ instead of $P_E$. Define $\tau:C(X)\ra
G$ by $\tau(\al)=P_\al/r_\al$, where $P_\al$ is the Hilbert
polynomial of $\al$, and $r_\al$ is the (positive) leading
coefficient of $P_\al$. Then $(\tau,G,\le)$ is a {\it permissible
stability condition\/} on $\coh(X)$ \cite[Ex.~4.16]{Joyc5}, called
{\it Gieseker stability}.

Gieseker stability is studied in \cite[\S 1.2]{HuLe}. Write
$\M_\ss^\al(\tau), \M_\st^\al(\tau)$ for the coarse moduli schemes
of $\tau$-(semi)stable sheaves $E$ with class $[E]=\al$ in $K(X)$.
By \cite[Th.~4.3.4]{HuLe}, $\M_\ss^\al(\tau)$ is a projective
$\C$-scheme whose $\C$-points correspond to S-equivalence classes of
Gieseker semistable sheaves in class $\al$, and $\M_\st^\al(\tau)$
is an open $\C$-subscheme whose $\C$-points correspond to
isomorphism classes of Gieseker stable sheaves in class~$\al$.
\label{ig3ex1}
\end{ex}

\begin{ex} In the situation of Example \ref{ig3ex1}, define
\begin{equation*}
M=\bigl\{p(t)=t^d+a_{d-1}t^{d-1}:d=0,1,2,3,\;\> a_{d-1}\in\Q,\;\>
a_{-1}=0\bigr\}\subset G
\end{equation*}
and restrict the total order $\le$ on $G$ to $M$. Define
$\mu:C(X)\ra M$ by $\mu(\al)=t^d+a_{d-1}t^{d-1}$ when
$\tau(\al)=P_\al/r_\al=t^d+a_{d-1}t^{d-1}+\cdots+a_0$, that is,
$\mu(\al)$ is the truncation of the polynomial $\tau(\al)$ in
Example \ref{ig3ex1} at its second term. Then as in
\cite[Ex.~4.17]{Joyc5}, $(\mu,M,\le)$ is a {\it permissible weak
stability condition\/} on $\coh(X)$. It is called $\mu$-{\it
stability}, and is studied in~\cite[\S 1.6]{HuLe}.
\label{ig3ex2}
\end{ex}

In \cite[\S 8]{Joyc5} we define interesting stack functions
$\bdss^\al(\tau),\bep^\al(\tau)$ in $\SFa(\fM)$.

\begin{dfn} Let $(\tau,T,\le)$ be a permissible weak stability
condition on $\coh(X)$. Define stack functions
$\bdss^\al(\tau)=\bde_{\fM_\ss^\al(\tau)}$ in $\SFa(\fM)$ for
$\al\in C(X)$. That is, $\bdss^\al(\tau)$ is the characteristic
function, in the sense of Definition \ref{ig2def5}, of the moduli
substack $\fM_\ss^\al(\tau)$ of $\tau$-semistable sheaves in $\fM$.
In \cite[Def.~8.1]{Joyc5} we define elements $\bep^\al(\tau)$ in
$\SFa(\fM)$ by
\ea
\bep^\al(\tau)&= \!\!\!\!\!\!\!
\sum_{\begin{subarray}{l}n\ge 1,\;\al_1,\ldots,\al_n\in C(X):\\
\al_1+\cdots+\al_n=\al,\; \tau(\al_i)=\tau(\al),\text{ all
$i$}\end{subarray}} \!\!\!\!\!\!
\frac{(-1)^{n-1}}{n}\,\,\bdss^{\al_1}(\tau)*\bdss^{\al_2}(\tau)
*\cdots*\bdss^{\al_n}(\tau),
\label{ig3eq6}
\intertext{where $*$ is the Ringel--Hall multiplication in
$\SFa(\fM)$. Then \cite[Th.~8.2]{Joyc5} proves} \bdss^\al(\tau)&=
\!\!\!\!\!\!\!
\sum_{\begin{subarray}{l}n\ge 1,\;\al_1,\ldots,\al_n\in C(X):\\
\al_1+\cdots+\al_n=\al,\; \tau(\al_i)=\tau(\al),\text{ all
$i$}\end{subarray}} \!\!\!
\frac{1}{n!}\,\,\bep^{\al_1}(\tau)*\bep^{\al_2}(\tau)*
\cdots*\bep^{\al_n}(\tau).
\label{ig3eq7}
\ea
There are only finitely many nonzero terms
in~\eq{ig3eq6}--\eq{ig3eq7}.
\label{ig3def5}
\end{dfn}

Equations \eq{ig3eq6} and \eq{ig3eq7} are inverse, so that knowing
the $\bep^\al(\tau)$ is equivalent to knowing the $\bdss^\al(\tau)$.
If $\fM_\ss^\al(\tau)=\fM_\st^\al(\tau)$ then
$\bep^\al(\tau)=\bdss^\al(\tau)$. The difference between
$\bep^\al(\tau)$ and $\bdss^\al(\tau)$ is that $\bep^\al(\tau)$
`counts' strictly semistable sheaves in a special, complicated way.
Here \cite[Th.~8.7]{Joyc5} is an important property of the
$\bep^\al(\tau)$, which does not hold for the $\bdss^\al(\tau)$. The
proof is highly nontrivial, using the full power of the
configurations formalism of~\cite{Joyc3,Joyc4,Joyc5,Joyc6}.

\begin{thm} $\bep^\al(\tau)$ lies in the Lie subalgebra
$\SFai(\fM)$ in~$\SFa(\fM)$.
\label{ig3thm2}
\end{thm}

In \cite[\S 6.6]{Joyc6} we define {\it invariants\/}
$J^\al(\tau)\in\Q$ for all $\al\in C(X)$ by
\e
\Psi\bigl(\bep^\al(\tau)\bigr)=J^\al(\tau)\la^\al.
\label{ig3eq8}
\e
This is valid by Theorem \ref{ig3thm2}. These $J^\al(\tau)$ are
rational numbers `counting' $\tau$-semistable sheaves $E$ in class
$\al$. When $\M_\ss^\al(\tau)=\M_\st^\al(\tau)$ we have
\e
J^\al(\tau)=\chi\bigl(\M_\st^\al(\tau)\bigr),
\label{ig3eq9}
\e
that is, $J^\al(\tau)$ is the Euler characteristic of the moduli
space $\M_\st^\al(\tau)$. In the notation of \S\ref{ig34}, this is
not weighted by the Behrend function $\nu_{\M_\st^\al(\tau)}$, and
so is not the Donaldson--Thomas invariant $DT^\al(\tau)$. As in
\cite[Ex.~6.9]{JoSo}, the $J^\al(\tau)$ are in general not unchanged
under deformations of~$X$.

\subsection{Changing stability conditions and algebra identities}
\label{ig33}

In \cite{Joyc6} we prove transformation laws for the
$\bdss^\al(\tau),\bep^\al(\tau)$ under change of stability
condition. These involve combinatorial coefficients
$S(*;\tau,\ti\tau)\in\Z$ and $U(*;\tau, \ti\tau)\in\Q$ defined in
\cite[\S 4.1]{Joyc6}.

\begin{dfn} Let $(\tau,T,\le),\!(\ti\tau,\ti T,\le)$ be weak
stability conditions on $\coh(X)$. Let $n\ge 1$ and
$\al_1,\ldots,\al_n\in C(X)$. If for all $i=1,\ldots,n-1$ we have
either
\begin{itemize}
\setlength{\itemsep}{0pt}
\setlength{\parsep}{0pt}
\item[(a)] $\tau(\al_i)\le\tau(\al_{i+1})$ and
$\ti\tau(\al_1+\cdots+\al_i)>\ti\tau(\al_{i+1}+\cdots+\al_n)$ or
\item[(b)] $\tau(\al_i)>\tau(\al_{i+1})$ and~
$\ti\tau(\al_1+\cdots+\al_i)\le\ti\tau(\al_{i+1}+\cdots+\al_n)$,
\end{itemize}
then define $S(\al_1,\ldots,\al_n;\tau,\ti\tau)=(-1)^r$, where $r$
is the number of $i=1,\ldots,n-1$ satisfying (a). Otherwise define
$S(\al_1,\ldots,\al_n;\tau,\ti\tau)=0$. Now define
\begin{align*}
&U(\al_1,\ldots,\al_n;\tau,\ti\tau)=\\
&\sum_{\begin{subarray}{l} \phantom{wiggle}\\
1\le l\le m\le n,\;\> 0=a_0<a_1<\cdots<a_m=n,\;\>
0=b_0<b_1<\cdots<b_l=m:\\
\text{Define $\be_1,\ldots,\be_m\in C(X)$ by
$\be_i=\al_{a_{i-1}+1}+\cdots+\al_{a_i}$.}\\
\text{Define $\ga_1,\ldots,\ga_l\in C(X)$ by
$\ga_i=\be_{b_{i-1}+1}+\cdots+\be_{b_i}$.}\\
\text{Then $\tau(\be_i)=\tau(\al_j)$, $i=1,\ldots,m$,
$a_{i-1}<j\le a_i$,}\\
\text{and $\ti\tau(\ga_i)=\ti\tau(\al_1+\cdots+\al_n)$,
$i=1,\ldots,l$}
\end{subarray}
\!\!\!\!\!\!\!\!\!\!\!\!\!\!\!\!\!\!\!\!\!\!\!\!\!\!\!\!\!\!\!\!\!
\!\!\!\!\!\!\!\!\!\!\!\!\!\!\!\!\!\!\!\!\!\!\!\!\!\!\!\!\!\!\!\!\!
\!\!\!\!\!\!\!\!\!\!\!\!\!\!\!\!\!\!\!\!}
\begin{aligned}[t]
\frac{(-1)^{l-1}}{l}\cdot\prod\nolimits_{i=1}^lS(\be_{b_{i-1}+1},
\be_{b_{i-1}+2},\ldots,\be_{b_i}; \tau,\ti\tau)&\\
\cdot\prod_{i=1}^m\frac{1}{(a_i-a_{i-1})!}&\,.
\end{aligned}
\end{align*}
\label{ig3def6}
\end{dfn}

Then in \cite[\S 5]{Joyc6} we derive wall-crossing formulae for the
$\bdss^\al(\tau),\bep^\al(\tau)$ under change of stability condition
from $(\tau,T,\le)$ to~$(\ti\tau,\ti T,\le)$:

\begin{thm} Let\/ $(\tau,T,\le), (\ti\tau,\ti T,\le)$ be
permissible weak stability conditions on $\coh(X)$. Then under some
mild extra conditions, for all\/ $\al\in C(X)$ we have
\begin{gather}
\begin{gathered}
\bdss^\al(\ti\tau)= \!\!\!\!\!\!\!
\sum_{\begin{subarray}{l}n\ge 1,\;\al_1,\ldots,\al_n\in
C(X):\\ \al_1+\cdots+\al_n=\al\end{subarray}} \!\!\!\!\!\!\!
\begin{aligned}[t]
S(\al_1,&\ldots,\al_n;\tau,\ti\tau)\cdot\\
&\bdss^{\al_1}(\tau)*\bdss^{\al_2}(\tau)*\cdots*
\bdss^{\al_n}(\tau),
\end{aligned}
\end{gathered}
\label{ig3eq10}\\
\begin{gathered}
\bep^\al(\ti\tau)= \!\!\!\!\!\!\!
\sum_{\begin{subarray}{l}n\ge 1,\;\al_1,\ldots,\al_n\in
C(X):\\ \al_1+\cdots+\al_n=\al\end{subarray}} \!\!\!\!\!\!\!
\begin{aligned}[t]
U(\al_1,&\ldots,\al_n;\tau,\ti\tau)\cdot\\
&\bep^{\al_1}(\tau)*\bep^{\al_2}(\tau)*\cdots* \bep^{\al_n}(\tau),
\end{aligned}
\end{gathered}
\label{ig3eq11}
\end{gather}
where there are only finitely many nonzero terms in
\eq{ig3eq10}--\eq{ig3eq11}.
\label{ig3thm3}
\end{thm}

The `mild extra conditions' in the theorem are required to ensure
that there are only finitely many nonzero terms in
\eq{ig3eq10}--\eq{ig3eq11}. In fact the author expects that this
always holds when $(\tau,T,\le),(\ti\tau,\ti T,\le)$ are of Gieseker
or $\mu$-stability type, but for irritating technical reasons has
not been able to prove this. As in \cite[\S 5.1]{Joyc6}, the author
can show that one can go between any two (weak) stability conditions
on $\coh(X)$ of Gieseker or $\mu$-stability type by finitely many
applications of Theorem \ref{ig3thm3}. In \cite[Th.~5.4]{Joyc6} we
prove:

\begin{thm} Equation \eq{ig3eq11} may be rewritten as an equation in
$\SFai(\fM)$ using the Lie bracket\/ $[\,,\,]$ on $\SFai(\fM),$
rather than as an equation in $\SFa(\fM)$ using the Ringel--Hall
product\/~$*$.
\label{ig3thm4}
\end{thm}

Therefore we may apply the Lie algebra morphism $\Psi$ of
\S\ref{ig31} to equation \eq{ig3eq11}. As \eq{ig3eq11} is not
expressed explicitly in terms of Lie brackets, it is helpful to
write this in the {\it universal enveloping algebra\/} $U(L(X))$.
This gives
\e
\begin{gathered}
J^\al(\ti\tau)\la^\al= \!\!\!\!\!\!\!
\sum_{\begin{subarray}{l}n\ge 1,\;\al_1,\ldots,\al_n\in
C(X):\\
\al_1+\cdots+\al_n=\al\end{subarray}\!\!\!\!\!\!\!\!\!\!\!\!\!\!\!
\!\!\!\!\!\!} \!\!\!\!\!\!\!\!
\begin{aligned}[t]
U(\al_1,\ldots,\al_n;\tau,\ti\tau)\,\cdot\, &\ts\prod_{i=1}^n
J^{\al_i}(\tau)\,\cdot\\
& \la^{\al_1}\star\la^{\al_2}\star\cdots\star\la^{\al_n},
\end{aligned}
\end{gathered}
\label{ig3eq12}
\e
where $\star$ is the product in~$U(L(X))$.

Now in \cite[\S 6.5]{Joyc4}, a basis is given for $U(L(X))$ in terms
of symbols $\la_{[I,\ka]}$, and multiplication $\star$ in $U(L(X))$
is written in terms of the $\la_{[I,\ka]}$ as a sum over graphs.
Here $I$ is a finite set, $\ka$ maps $I\ra C(X)$, and when
$\md{I}=1$, so that $I=\{i\}$, we have $\la_{[I,\ka]}=\la^{\ka(i)}$.
Then \cite[eq.~(127)]{Joyc6} gives an expression for
$\la^{\al_1}\star\cdots\star \la^{\al_n}$ in $U(L(X))$, in terms of
sums over {\it directed graphs\/} ({\it digraphs\/}):
\begin{gather}
\la^{\al_1}\star\cdots\star\la^{\al_n}=\text{ terms in
$\la_{[I,\ka]}$, $\md{I}>1$, }
\label{ig3eq13}\\
+\raisebox{-6pt}{\begin{Large}$\displaystyle\biggl[$\end{Large}}
\frac{1}{2^{n-1}}\!\!\!\!\!
\sum_{\substack{\text{connected, simply-connected digraphs
$\Ga$:}\\
\text{vertices $\{1,\ldots,n\}$, edge $\mathop{\bu} \limits^{\sst
i}\ra\mathop{\bu}\limits^{\sst j}$ implies $i<j$}}} \,\,\,
\prod_{\substack{\text{edges}\\
\text{$\mathop{\bu}\limits^{\sst i}\ra\mathop{\bu}\limits^{\sst
j}$}\\ \text{in $\Ga$}}}\bar\chi(\al_i,\al_j)
\raisebox{-6pt}{\begin{Large}$\displaystyle\biggr]$\end{Large}}
\la^{\al_1+\cdots+\al_n}. \nonumber
\end{gather}

Substitute \eq{ig3eq13} into \eq{ig3eq12}. The terms in
$\la_{[I,\ka]}$ for $\md{I}>1$ all cancel, as \eq{ig3eq12} lies in
$L(X)\subset U(L(X))$. So equating coefficients of $\la^\al$ yields
\e
\begin{gathered}
J^\al(\ti\tau)=\!\!\!\!\!\!
\sum_{\begin{subarray}{l}n\ge 1,\;\al_1,\ldots,\al_n\in
C(X):\\ \al_1+\cdots+\al_n=\al\end{subarray}}\,\,\,\,
\sum_{\begin{subarray}{l}\text{connected, simply-connected digraphs
$\Ga$:}\\ \text{vertices $\{1,\ldots,n\}$, edge $\mathop{\bu}
\limits^{\sst i}\ra\mathop{\bu}\limits^{\sst j}$ implies
$i<j$}\end{subarray}}\\ \frac{1}{2^{n-1}}\,
U(\al_1,\ldots,\al_n;\tau,\ti\tau) \!\!\!\!\! \prod_{\text{edges
$\mathop{\bu}\limits^{\sst i}\ra\mathop{\bu}\limits^{\sst j}$ in
$\Ga$}}\!\!\!\!\! \bar\chi(\al_i,\al_j)
\prod_{i=1}^nJ^{\al_i}(\tau).
\end{gathered}
\label{ig3eq14}
\e

Following \cite[Def.~6.27]{Joyc6}, we define combinatorial
coefficients~$V(I,\Ga,\ka;\tau,\ti\tau)$:

\begin{dfn} In the situation above, let $\Ga$ be a connected,
simply-connected digraph with finite vertex set $I$, where
$\md{I}=n$, and $\ka:I\ra C(X)$ be a map. Define
$V(I,\Ga,\ka;\tau,\ti\tau)\in\Q$ by
\e
V(I,\Ga,\ka;\tau,\ti\tau)=\frac{1}{2^{n-1}n!}
\!\!\sum_{\substack{\text{orderings $i_1,\ldots,i_n$ of $I$:}\\
\text{edge $\mathop{\bu} \limits^{\sst
i_a}\ra\mathop{\bu}\limits^{\sst i_b}$ in $\Ga$ implies $a<b$}}
\!\!\!\!\!\!\!\!\!\!\!\!\!\!\!\!\!\!\!\!\!\!\!\!\!\!\!\!\!\!\!\!\!\!\!
} \!\!\!\!\! U(\ka(i_1),\ka(i_2),\ldots,\ka(i_n);\tau,\ti\tau).
\label{ig3eq15}
\e
\label{ig3def7}
\end{dfn}

Then as in \cite[Th.~6.28]{Joyc6}, using \eq{ig3eq15} to rewrite
\eq{ig3eq14} yields a transformation law for the $J^\al(\tau)$ under
change of stability condition:
\e
\begin{gathered}
J^\al(\ti\tau)\!=\!\!\!\!
\sum_{\substack{\text{iso.}\\ \text{classes}\\
\text{of finite}\\ \text{sets $I$}}}\,\,
\sum_{\substack{\ka:I\ra C(X):\\ \sum_{i\in I}\ka(i)=\al}}\,\,
\sum_{\begin{subarray}{l} \text{connected,}\\
\text{simply-connected}\\ \text{digraphs $\Ga$,}\\
\text{vertices $I$}\end{subarray}\!\!\!\!\!\!\!\!\!\!\!\!\!\!\!\!\!
\!\!\!\!\!\!\!\!\!\!} V(I,\Ga,\ka;\tau,\ti\tau)
\begin{aligned}[t]
&\cdot\prod\limits_{\text{edges \smash{$\mathop{\bu}\limits^{\sst
i}\ra\mathop{\bu}\limits^{\sst j}$} in
$\Ga$}\!\!\!\!\!\!\!\!\!\!\!\!\!\!\!\!\!\!\!\!\!\!\!
\!\!\!\!\!\!\!\!\!\!\!\!\!\!} \bar\chi(\ka(i),\ka(j))\\
&\cdot\prod\nolimits_{i\in I}J^{\ka(i)}(\tau).
\end{aligned}
\end{gathered}
\label{ig3eq16}
\e

\subsection{Donaldson--Thomas invariants of Calabi--Yau 3-folds}
\label{ig34}

{\it Donaldson--Thomas invariants\/} $DT^\al(\tau)$ were defined by
Richard Thomas \cite{Thom}, following a proposal of Donaldson and
Thomas~\cite[\S 3]{DoTh}.

\begin{dfn} Let $X$ be a Calabi--Yau 3-fold. Fix a very ample line bundle
$\O_X(1)$ on $X$, and let $(\tau,G,\le)$ be Gieseker stability on
$\coh(X)$ w.r.t.\ $\O_X(1)$, as in Example \ref{ig3ex1}. For $\al\in
K(X)$, write $\M_\ss^\al(\tau),\M_\st^\al(\tau)$ for the coarse
moduli schemes of $\tau$-(semi)stable sheaves $E$ with class
$[E]=\al$. Then $\M_\ss^\al(\tau)$ is a projective $\C$-scheme, and
$\M_\st^\al(\tau)$ an open subscheme.

Thomas \cite{Thom} constructs a symmetric obstruction theory on
$\M_\st^\al(\tau)$. Suppose that $\M_\ss^\al(\tau)=
\M_\st^\al(\tau)$. Then $\M_\st^\al(\tau)$ is proper, so using the
obstruction theory Behrend and Fantechi \cite{BeFa1} define a
virtual class $[\M_\st^\al(\tau)]^\vir\in A_0(\M_\st^\al(\tau))$.
The {\it Donaldson--Thomas invariant\/} \cite{Thom} is defined to be
\e
DT^\al(\tau)=\ts\int_{[\M_\st^\al(\tau)]^\vir}1.
\label{ig3eq17}
\e
Note that $DT^\al(\tau)$ {\it is defined only when\/ $\M_\ss^\al
(\tau)=\M_\st^\al(\tau)$, that is, there are no strictly semistable
sheaves\/ $E$ in class\/} $\al$. One of our main goals is to extend
the definition to all $\al\in K(X)$. Thomas' main result \cite[\S
3]{Thom} is that
\label{ig3def8}
\end{dfn}

\begin{thm} $DT^\al(\tau)$ is unchanged by continuous
deformations of the underlying Calabi--Yau $3$-fold~$X$.
\label{ig3thm5}
\end{thm}

An important advance in Donaldson--Thomas theory was made by Behrend
\cite{Behr}, who found a way to rewrite the definition \eq{ig3eq17}
of Donaldson--Thomas invariants as a weighted Euler characteristic.
Let $\fF$ be an Artin $\C$-stack, locally of finite type. Then $\fF$
has a unique {\it Behrend function\/} $\nu_\fF:\fF(\C)\ra\Z,$ a
$\Z$-valued locally constructible function on $\fF$. The definition,
which we do not give, can be found in \cite[\S 1]{Behr} when $\fF$
is a finite type $\C$-scheme, and in \cite[\S 4.1]{JoSo} in the
general case. Here are some important properties of Behrend
functions from~\cite{Behr,JoSo}.

\begin{thm} Let\/ $\fF,\fG$ be Artin $\C$-stacks locally of
finite type. Then:
\begin{itemize}
\setlength{\itemsep}{0pt}
\setlength{\parsep}{0pt}
\item[{\rm(i)}] If\/ $\fF$ is a smooth of dimension\/ $n$
then~$\nu_\fF\equiv(-1)^n$.
\item[{\rm(ii)}] If\/ $\vp:\fF\!\ra\! \fG$ is smooth with
relative dimension $n$
then\/~$\nu_\fF\!\equiv\!(-1)^n\vp^*(\nu_\fG)$.
\item[{\rm(iii)}] $\nu_{\fF\t \fG}=\nu_\fF\boxdot\nu_\fG$ in
$\LCF(\fF\t \fG),$ where $(\nu_\fF\boxdot\nu_\fG)(x,y)=
\nu_\fF(x)\nu_\fG(y)$.
 \item[{\rm(iv)}] Suppose $\M$ is a proper $\C$-scheme and has a
symmetric obstruction theory, and\/ $[\M]^\vir\in A_0(\M)$ is the
corresponding virtual class from\/ {\rm\cite{BeFa1}}. Then
\begin{equation*}
\ts\int_{[\M]^\vir}1=\chi(\M,\nu_\M)\in\Z,
\end{equation*}
where $\chi(\M,\nu_\M)=\int_{\M(\C)}\nu_\M\d\chi$ is the
\begin{bfseries}weighted Euler characteristic\end{bfseries} of\/
$\M,$ weighted by the constructible function $\nu_\M$. In
particular, $\int_{[\M]^\vir}1$ depends only on the\/ $\C$-scheme
structure of\/ $\M,$ not on the choice of symmetric obstruction
theory.
\item[{\rm(v)}] Let\/ $\M$ be a $\C$-scheme, let\/ $x\in
\M(\C),$ and suppose there exist a complex manifold\/ $U,$ a
holomorphic function\/ $f:U\ra\C,$ and a point\/
$u\in\Crit(f)\subseteq U$ such that locally in the analytic
topology, $\M(\C)$ near $x$ is isomorphic as a complex analytic
space to\/ $\Crit(f)$ near\/ $u$. Then
\begin{equation*}
\nu_\M(x)=(-1)^{\dim U}\bigl(1-\chi(MF_f(u))\bigr),
\end{equation*}
where $\chi(MF_f(u))$ is the Euler characteristic of the Milnor
fibre $MF_f(u)$.
\end{itemize}
\label{ig3thm6}
\end{thm}

Here the {\it Milnor fibre\/} in (v) is defined as follows:

\begin{dfn} Let $U$ be a complex analytic space, locally of finite
type, $f:U\ra\C$ a holomorphic function, and $u\in U$. Let
$d(\,,\,)$ be a metric on $U$ near $u$ induced by a local embedding
of $U$ in some $\C^N$. For $u\in U$ and $\de,\ep>0$, consider the
holomorphic map
\begin{equation*}
\Phi_{f,u}:\bigl\{v\in U:d(u,v)\!<\!\de,\;
0\!<\!\md{f(v)\!-\!f(u)}\!<\!\ep\bigr\}\longra
\bigl\{z\in\C:0\!<\!\md{z}\!<\!\ep\bigr\}
\end{equation*}
given by $\Phi_{f,u}(v)=f(v)-f(u)$. Then $\Phi_{f,u}$ is a smooth
locally trivial fibration provided $0<\ep\ll\de\ll 1$. The {\it
Milnor fibre\/} $MF_f(u)$ is the fibre of $\Phi_{f,u}$. It is
independent of the choice of~$0<\ep\ll\de\ll 1$.
\label{ig3def9}
\end{dfn}

Theorem \ref{ig3thm6}(iv) implies that $DT^\al(\tau)$ in
\eq{ig3eq17} is given by
\e
DT^\al(\tau)=\chi\bigl(\M_\st^\al(\tau),\nu_{\M_\st^\al(\tau)}\bigr).
\label{ig3eq18}
\e
This is similar to the expression \eq{ig3eq9} for $J^\al(\tau)$ when
$\M_\ss^\al(\tau)= \M_\st^\al(\tau)$. There is a big difference
between the two equations \eq{ig3eq17} and \eq{ig3eq18} defining
Donaldson--Thomas invariants. Equation \eq{ig3eq17} is non-local,
and non-motivic, and makes sense only if $\M_\st^\al(\tau)$ is a
proper $\C$-scheme. But \eq{ig3eq18} is local, and (in a sense)
motivic, and makes sense for arbitrary finite type $\C$-schemes
$\M_\st^\al(\tau)$. It is tempting to take \eq{ig3eq18} to be the
definition of Donaldson--Thomas invariants even when
$\M_\ss^\al(\tau)\ne\M_\st^\al(\tau)$, but in \cite[\S 6.5]{JoSo} we
show that this is not a good idea, as then $DT^\al(\tau)$ would not
be unchanged under deformations of~$X$.

Equation \eq{ig3eq18} was the inspiration for \cite{JoSo}. It shows
that Donaldson--Thomas invariants $DT^\al(\tau)$ can be written as
{\it motivic\/} invariants, like those studied in
\cite{Joyc3,Joyc4,Joyc5,Joyc6,Joyc7}, and suggests extending the
results of \cite{Joyc3,Joyc4,Joyc5,Joyc6,Joyc7} to Donaldson--Thomas
invariants by including Behrend functions as weights.

\section{Generalized Donaldson--Thomas invariants}
\label{ig4}

We now summarize \cite[\S 5--\S 6]{JoSo}. All this section is joint
work with Yinan Song. Let $X$ be a Calabi--Yau 3-fold over $\C$, and
$\O_X(1)$ a very ample line bundle over $X$. We now assume that
$H^1(\O_X)=0$, which was not needed in \S\ref{ig3}. We use the
notation of \S\ref{ig3}, with $\fM$ the moduli stack of coherent
sheaves on $X$, and so on.

\subsection{Local description of the moduli of coherent sheaves}
\label{ig41}

In \cite[Th.~5.5]{JoSo} we give a local characterization of an atlas
for the moduli stack $\fM$ as the critical points of a holomorphic
function on a complex manifold.

\begin{thm} Let\/ $X$ be a Calabi--Yau $3$-fold over\/ $\C,$ and\/
$\fM$ the moduli stack of coherent sheaves on\/ $X$. Suppose\/ $E$
is a coherent sheaf on\/ $X,$ so that\/ $[E]\in\fM(\C)$. Let\/ $G$
be a maximal compact subgroup in $\Aut(E),$ and\/ $G^{\sst\C}$ its
complexification. Then\/ $G^{\sst\C}$ is an algebraic $\C$-subgroup
of\/ $\Aut(E),$ a maximal reductive subgroup, and\/
$G^{\sst\C}=\Aut(E)$ if and only if\/ $\Aut(E)$ is reductive.

There exists a quasiprojective $\C$-scheme $S,$ an action of\/
$G^{\sst\C}$ on $S,$ a point\/ $s\in S(\C)$ fixed by $G^{\sst\C},$
and a $1$-morphism of Artin $\C$-stacks $\Phi:[S/G^{\sst\C}]\ra\fM,$
which is smooth of relative dimension $\dim\Aut(E)-\dim G^{\sst\C},$
where $[S/G^{\sst\C}]$ is the quotient stack, such that\/
$\Phi(s\,G^{\sst\C})=[E],$ the induced morphism on stabilizer groups
$\Phi_*:\Iso_{[S/G^{\sst\C}]}(s\,G^{\sst\C})\ra\Iso_{\fM}([E])$ is
the natural morphism $G^{\sst\C}\hookra\Aut(E)\cong\Iso_{\fM}([E]),$
and\/ $\d\Phi\vert_{s\,G^{\sst\C}}:T_sS\cong T_{s\,G^{\sst\C}}
[S/G^{\sst\C}]\ra T_{[E]}\fM\cong \Ext^1(E,E)$ is an isomorphism.
Furthermore, $S$ parametrizes a formally versal family $(S,{\cal
D})$ of coherent sheaves on $X,$ equivariant under the action of\/
$G^{\sst\C}$ on $S,$ with fibre\/ ${\cal D}_s\cong E$ at\/ $s$. If\/
$\Aut(E)$ is reductive then $\Phi$ is \'etale.

Write $S_\an$ for the complex analytic space underlying the
$\C$-scheme $S$. Then there exists an open neighbourhood\/ $U$ of\/
$0$ in\/ $\Ext^1(E,E)$ in the analytic topology, a holomorphic
function $f:U\ra\C$ with\/ $f(0)=\d f\vert_0=0,$ an open
neighbourhood\/ $V$ of\/ $s$ in $S_\an,$ and an isomorphism of
complex analytic spaces $\Xi:\Crit(f)\ra V,$ such that\/ $\Xi(0)=s$
and\/ $\d\Xi\vert_0:T_0\Crit(f)\ra T_sV$ is the inverse of\/
$\d\Phi\vert_{s\,G^{\sst\C}}:T_sS\ra\Ext^1(E,E)$. Moreover we can
choose $U,f,V$ to be $G^{\sst\C}$-invariant, and\/ $\Xi$ to be
$G^{\sst\C}$-equivariant.
\label{ig4thm1}
\end{thm}

The proof of Theorem \ref{ig4thm1} comes in two parts. First we show
in \cite[\S 8]{JoSo} that $\fM$ near $[E]$ is locally isomorphic, as
an Artin $\C$-stack, to the moduli stack $\fVect$ of algebraic
vector bundles on $X$ near $[E']$ for some vector bundle $E'\ra X$.
The proof uses algebraic geometry, and is valid for $X$ an
Calabi--Yau $m$-fold for any $m>0$ over any algebraically closed
field $\K$. The local morphism $\fM\ra\fVect$ is the composition of
shifts and $m$ Seidel--Thomas twists by $\O_X(-n)$ for~$n\gg 0$.

Thus, it is enough to prove Theorem \ref{ig4thm1} with $\fVect$ in
place of $\fM$. We do this in \cite[\S 9]{JoSo} using gauge theory
on vector bundles over $X$, motivated by an idea of Donaldson and
Thomas \cite[\S 3]{DoTh}, \cite[\S 2]{Thom}, and results of Miyajima
\cite{Miya}. Let $E\ra X$ be a fixed complex (not holomorphic)
vector bundle over $X$. Write $\sA$ for the infinite-dimensional
affine space of smooth semiconnections ($\db$-operators) on $E$, and
$\sG$ for the infinite-dimensional Lie group of smooth gauge
transformations of $E$. Then $\sG$ acts on $\sA$, and $\sB=\sA/\sG$
is the space of gauge-equivalence classes of semiconnections on~$E$.

We fix $\db_E$ in $\sA$ coming from a holomorphic vector bundle
structure on $E$. Then points in $\sA$ are of the form $\db_E+A$ for
$A\in C^\iy\bigl(\End(E)\ot_\C \La^{0,1}T^*X\bigr)$, and $\db_E+A$
makes $E$ into a holomorphic vector bundle if $F_A^{0,2}=\db_EA+A\w
A$ is zero in $\smash{C^\iy\bigl(\End(E)\ot_\C
\La^{0,2}T^*X\bigr)}$. Thus, the moduli space (stack) of holomorphic
vector bundle structures on $E$ is isomorphic to $\{\db_E+A\in\sA:
F_A^{0,2}=0\}/\sG$. Thomas observes that when $X$ is a Calabi--Yau
3-fold, there is a natural holomorphic function $CS:\sA\ra\C$ called
the {\it holomorphic Chern--Simons functional}, invariant under
$\sG$ up to addition of constants, such that
$\{\db_E+A\in\sA:F_A^{0,2}=0\}$ is the critical locus of $CS$. Thus,
$\fVect$ is (informally) locally the critical points of a
holomorphic function $CS$ on an infinite-dimensional complex stack
$\sB=\sA/\sG$. To prove Theorem \ref{ig4thm1} we show that we can
find a finite-dimensional complex submanifold $U$ in $\sA$ and a
finite-dimensional complex Lie subgroup $G^{\sst\C}$ in $\sG$
preserving $U$ such that the theorem holds with~$f=CS\vert_U$.

In \cite[Th.~5.11]{JoSo} we prove identities on the {\it Behrend
function\/} of $\fM$, as in~\S\ref{ig34}.

\begin{thm} Let\/ $X$ be a Calabi--Yau $3$-fold over\/ $\C,$ and\/
$\fM$ the moduli stack of coherent sheaves on\/ $X$. The
\begin{bfseries}Behrend function\end{bfseries}\/ $\nu_{\fM}:
\fM(\C)\ra\Z$ is a natural locally constructible function on $\fM$.
For all\/ $E_1,E_2\in\coh(X),$ it satisfies:
\begin{gather}
\nu_{\fM}(E_1\op E_2)=(-1)^{\bar\chi([E_1],[E_2])}
\nu_{\fM}(E_1)\nu_{\fM}(E_2),
\label{ig4eq1}\\
\begin{split}
\int_{\begin{subarray}{l}[\la]\in\mathbb{P}(\Ext^1(E_2,E_1)):\\
\la\; \Leftrightarrow\; 0\ra E_1\ra F\ra E_2\ra
0\end{subarray}}\!\!\!\!\!\!\nu_{\fM}(F)\,\d\chi -
\int_{\begin{subarray}{l}[\la']\in\mathbb{P}(\Ext^1(E_1,E_2)):\\
\la'\; \Leftrightarrow\; 0\ra E_2\ra F'\ra E_1\ra
0\end{subarray}}\!\!\!\!\!\!\nu_{\fM}(F')\,\d\chi \\
=\bigl(\dim\Ext^1(E_2,E_1)-\dim\Ext^1(E_1,E_2)\bigr)
\nu_{\fM}(E_1\op E_2).
\end{split}
\label{ig4eq2}
\end{gather}
Here\/ $\bar\chi([E_1],[E_2])$ in \eq{ig4eq1} is defined in
{\rm\eq{ig3eq1},} and in \eq{ig4eq2} the correspondence between\/
$[\la]\in\mathbb{P}(\Ext^1(E_2,E_1))$ and\/ $F\in\coh(X)$ is that\/
$[\la]\in\mathbb{P}(\Ext^1(E_2,E_1))$ lifts to some\/
$0\ne\la\in\Ext^1(E_2,E_1),$ which corresponds to a short exact
sequence\/ $0\ra E_1\ra F\ra E_2\ra 0$ in\/ $\coh(X)$ in the usual
way. The function $[\la]\mapsto\nu_{\fM}(F)$ is a constructible
function\/ $\mathbb{P}(\Ext^1(E_2,E_1))\ra\Z,$ and the integrals in
\eq{ig4eq2} are integrals of constructible functions using the Euler
characteristic as measure.
\label{ig4thm2}
\end{thm}

We prove Theorem \ref{ig4thm2} using Theorem \ref{ig4thm1} and the
Milnor fibre description of Behrend functions from Theorem
\ref{ig3thm6}(v). We apply Theorem \ref{ig4thm1} to $E=E_1\op E_2$,
and we take the maximal compact subgroup $G$ of $\Aut(E)$ to contain
the subgroup $\bigl\{\id_{E_1}+\la\id_{E_2}: \la\in\U(1)\bigr\}$, so
that $G^{\sst\C}$ contains $\bigl\{\id_{E_1}+\la\id_{E_2}:\la\in
\bG_m\bigr\}$. Equations \eq{ig4eq1} and \eq{ig4eq2} are proved by a
kind of localization using this $\bG_m$-action on~$\Ext^1(E_1\op
E_2,E_1\op E_2)$.

Note that Theorem \ref{ig4thm2} makes sense as a statement in
algebraic geometry, for Calabi--Yau 3-folds over an algebraically
closed field $\K$ of characteristic zero, and the author expects it
to be true in this generality. However, our proof of Theorem
\ref{ig4thm2} uses gauge theory, and transcendental complex analytic
geometry methods, and is valid only over~$\K=\C$.

\subsection{A Lie algebra morphism $\ti\Psi:\SFai(\fM)\ra\ti L(X),$
and \\ generalized Donaldson--Thomas invariants
$\bar{DT}{}^\al(\tau)$}
\label{ig42}

In \S\ref{ig31} we defined an explicit Lie algebra $L(X)$ and Lie
algebra morphisms $\Psi:\SFai(\fM)\ra L(X)$ and
$\Psi^{\chi,\Q}:\oSFai (\fM,\chi,\Q)\ra L(X)$. We now define
modified versions $\ti L(X),\ti\Psi,\ti\Psi^{\chi,\Q}$, with
$\ti\Psi,\ti\Psi^{\chi,\Q}$ weighted by the Behrend function
$\nu_\fM$ of $\fM$. We continue to use the notation
of~\S\ref{ig2}--\S\ref{ig3}.

\begin{dfn} Define a Lie algebra $\ti L(X)$ to be the $\Q$-vector
space with basis of symbols $\ti \la^\al$ for $\al\in K(X)$, with
Lie bracket
\e
[\ti\la^\al,\ti \la^\be]=(-1)^{\bar\chi(\al,\be)}\bar\chi(\al,\be)
\ti\la^{\al+\be},
\label{ig4eq3}
\e
which is \eq{ig3eq3} with a sign change. As $\bar\chi$ is
antisymmetric, \eq{ig4eq3} satisfies the Jacobi identity, and makes
$\ti L(X)$ into an infinite-dimensional Lie algebra over~$\Q$.

Define a $\Q$-linear map $\ti\Psi^{\chi,\Q}:\oSFai
(\fM,\chi,\Q)\ra\ti L(X)$ by
\begin{equation*}
\ti\Psi^{\chi,\Q}(f)=\ts\sum_{\al\in K(X)}\ga^\al\ti \la^{\al},
\end{equation*}
as in \eq{ig3eq4}, where $\ga^\al\in\Q$ is defined as follows. Write
$f\vert_{\fM^\al}$ in terms of $\de_i,U_i,\rho_i$ as in \eq{ig3eq5},
and set
\e
\ga^\al=\ts\sum_{i=1}^n\de_i\chi\bigl(U_i,\rho_i^*(\nu_\fM)\bigr),
\label{ig4eq4}
\e
where $\rho_i^*(\nu_\fM)$ is the pullback of the Behrend function
$\nu_\fM$ to a constructible function on $U_i\t[\Spec\C/\bG_m]$, or
equivalently on $U_i$, and $\chi\bigl(U_i,\rho_i^*(\nu_\fM)\bigr)$
is the Euler characteristic of $U_i$ weighted by
$\rho_i^*(\nu_\fM)$. One can show that the map from \eq{ig3eq5} to
\eq{ig4eq4} is compatible with the relations in
$\oSFai(\fM^\al,\chi,\Q)$, and so $\ti\Psi^{\chi,\Q}$ is
well-defined. Define $\ti\Psi:\SFai(\fM)\ra\ti L(X)$
by~$\ti\Psi=\ti\Psi^{\chi,\Q}\ci\bar\Pi^{\chi,\Q}_\fM$.
\label{ig4def1}
\end{dfn}

The reason for the sign change between \eq{ig3eq3} and \eq{ig4eq3}
is the signs involved in Behrend functions, in particular, the
$(-1)^n$ in Theorem \ref{ig3thm6}(ii), which is responsible for the
factor $(-1)^{\bar\chi([E_1],[E_2])}$ in \eq{ig4eq1}. Here
\cite[Th.~5.14]{JoSo} is the analogue of Theorem~\ref{ig3thm1}.

\begin{thm} $\ti\Psi:\SFai(\fM)\ra\ti L(X)$ and\/
$\ti\Psi^{\chi,\Q}:\oSFai(\fM,\chi,\Q)\ra\ti L(X)$ are Lie algebra
morphisms.
\label{ig4thm3}
\end{thm}

We can now define generalized Donaldson--Thomas invariants.

\begin{dfn} Let $X$ be a projective Calabi--Yau 3-fold over $\C$,
let $\O_X(1)$ be a very ample line bundle on $X$, and let
$(\tau,G,\le)$ be Gieseker stability and $(\mu,M,\le)$ be
$\mu$-stability on $\coh(X)$ w.r.t.\ $\O_X(1)$, as in Examples
\ref{ig3ex1} and \ref{ig3ex2}. As in \eq{ig3eq8}, define {\it
generalized Donaldson--Thomas invariants\/}
$\bar{DT}{}^\al(\tau)\in\Q$ and $\bar{DT}{}^\al(\mu)\in\Q$ for all
$\al\in C(X)$ by
\e
\ti\Psi\bigl(\bep^\al(\tau)\bigr)=-\bar{DT}{}^\al(\tau)\ti
\la^\al\qquad\text{and}\qquad
\ti\Psi\bigl(\bep^\al(\mu)\bigr)=-\bar{DT}{}^\al(\mu)\ti \la^\al.
\label{ig4eq5}
\e
Here $\bep^\al(\tau),\bep^\al(\mu)$ are defined in \eq{ig3eq6}, and
lie in $\SFai(\fM)$ by Theorem \ref{ig3thm2}, so
$\bar{DT}{}^\al(\tau),\bar{DT}{}^\al(\mu)$ are well-defined. In
\cite[Prop.~5.17]{JoSo} we show that if $\M_\ss^\al(\tau)=
\M_\st^\al(\tau)$ then $\bar{DT}{}^\al(\tau)=DT^\al(\tau)$. That is,
our new generalized Donaldson--Thomas invariants
$\bar{DT}{}^\al(\tau)$ are equal to the original Donaldson--Thomas
invariants $DT^\al(\tau)$ of \cite{Thom} whenever the $DT^\al(\tau)$
are defined.
\label{ig4def2}
\end{dfn}

We can now repeat the argument of \S\ref{ig33} to deduce
transformation laws for generalized Donaldson--Thomas invariants
under change of stability condition. In the situation of Theorem
\ref{ig3thm3}, equation \eq{ig3eq11} is an identity in the Lie
algebra $\SFai(\fM)$, so we can apply the Lie algebra morphism
$\ti\Psi$ to transform \eq{ig3eq11} into an identity in the Lie
algebra $\ti L(X)$, and use \eq{ig4eq5} to write this in terms of
generalized Donaldson--Thomas invariants. As for \eq{ig3eq12}, this
gives an equation in the universal enveloping algebra~$U(\ti L(X))$:
\begin{equation*}
\bar{DT}{}^\al(\ti\tau)\ti \la^\al= \!\!\!\!\!\!\!
\sum_{\begin{subarray}{l}n\ge 1,\;\al_1,\ldots,\al_n\in
C(X):\\ \al_1+\cdots+\al_n=\al\end{subarray}
\!\!\!\!\!\!\!\!\!\!\!\!\!\!\!\!\!\!\!\!\!} \!\!\!\!\!\!\!\!\!
\begin{aligned}[t]
U(\al_1,\ldots,\al_n;\tau,\ti\tau)\,\cdot\, &\ts
(-1)^{n-1}\prod_{i=1}^n\bar{DT}{}^{\al_i}(\tau)\cdot\\
& \ti \la^{\al_1}\star \ti \la^{\al_2}\star\cdots\star\ti
\la^{\al_n}.
\end{aligned}
\end{equation*}
Following the proof of \eq{ig3eq16} in \S\ref{ig33} with sign
changes, in \cite[Th.~5.18]{JoSo} we obtain:

\begin{thm} In the situation of Theorem\/
{\rm\ref{ig3thm3},} for all\/ $\al\in C(X)$ we have
\e
\begin{aligned}
&\bar{DT}{}^\al(\ti\tau)=\\
&\!\!\!\sum_{\substack{\text{iso.}\\ \text{classes}\\
\text{of finite}\\ \text{sets $I$}}}\,
\sum_{\substack{\ka:I\ra C(X):\\ \sum\limits_{i\in I}\ka(i)=\al}}\,\,
\sum_{\begin{subarray}{l} \text{connected,}\\
\text{simply-}\\ \text{connected}\\ \text{digraphs $\Ga$,}\\
\text{vertices $I$}\end{subarray}}
\begin{aligned}[t]
(-1)^{\md{I}-1} V(I,\Ga,\ka;\tau,\ti\tau)\cdot
\prod\nolimits_{i\in I} \bar{DT}{}^{\ka(i)}(\tau)&\\
\cdot (-1)^{\frac{1}{2}\sum_{i,j\in
I}\md{\bar\chi(\ka(i),\ka(j))}}\cdot\! \prod\limits_{\text{edges
\smash{$\mathop{\bu}\limits^{\sst i}\ra\mathop{\bu}\limits^{\sst
j}$} in $\Ga$}\!\!\!\!\!\!\!\!\!\!\!\!\!\!\!\!\!\!\!\!\!\!\!
\!\!\!\!\!\!\!\!} \bar\chi(\ka(i),\ka(j))&,
\end{aligned}
\end{aligned}
\label{ig4eq6}
\e
with only finitely many nonzero terms.
\label{ig4thm4}
\end{thm}

The discussion after Theorem \ref{ig3thm3}
implies~\cite[Cor.~5.19]{JoSo}:

\begin{cor} Let\/ $(\tau,T,\le),(\ti\tau,\ti T,\le)$ be two
permissible weak stability conditions on $\coh(X)$ of Gieseker or
$\mu$-stability type, as in Examples {\rm\ref{ig3ex1}} and\/
{\rm\ref{ig3ex2}}. Then the $\bar{DT}{}^\al(\tau)$ for all\/ $\al\in
C(X)$ completely determine the $\bar{DT}{}^\al(\ti\tau)$ for all\/
$\al\in C(X),$ and vice versa, through finitely many applications
of\/~\eq{ig4eq6}.
\label{ig4cor1}
\end{cor}

\subsection{Invariants $PI^{\al,n}(\tau')$ counting stable pairs,
and \\ deformation-invariance of the $\bar{DT}{}^\al(\tau)$}
\label{ig43}

We wish to prove that our invariants $\bar{DT}{}^\al(\tau)$ are
unchanged under deformations of $X$. We do this indirectly: we first
define another family of auxiliary invariants $PI^{\al,n}(\tau')$
counting {\it stable pairs\/} on $X$, and show that
$PI^{\al,n}(\tau')$ are unchanged under deformations of $X$. Then we
prove an identity \eq{ig4eq9} expressing $PI^{\al,n}(\tau')$ in
terms of the $\bar{DT}{}^\be(\tau)$, and use it to show
$\bar{DT}{}^\al(\tau)$ is deformation-invariant. This approach was
inspired by Pandharipande and Thomas \cite{PaTh}, who use invariants
counting pairs to study curve counting in Calabi--Yau 3-folds.

\begin{dfn} Let $X$ be a Calabi--Yau 3-fold over $\C$, with
$H^1(\O_X)=0$. Choose a very ample line bundle $\O_X(1)$ on $X$, and
write $(\tau,G,\le)$ for Gieseker stability w.r.t.\ $\O_X(1)$, as in
Example~\ref{ig3ex1}.

Fix $n\gg 0$ in $\Z$. A {\it pair\/} is a nonzero morphism of
sheaves $s:\O_{X}(-n)\ra E$, where $E$ is a nonzero sheaf. A {\it
morphism\/} between two pairs $s:\O_{X}(-n)\ra E$ and
$t:\O_{X}(-n)\ra F$ is a morphism of $\O_X$-modules $f:E\ra F$, with
$f\ci s=t$. A pair $s:\O_X(-n)\ra E$ is called {\it stable\/} if:
\begin{itemize}
\setlength{\itemsep}{0pt}
\setlength{\parsep}{0pt}
\item[(i)] $\tau([E'])\le\tau([E])$ for all subsheaves $E'$ of
$E$ with $0\neq E'\neq E$; and
\item[(ii)] If also $s$ factors through $E'$,
then~$\tau([E'])<\tau([E])$.
\end{itemize}
Note that (i) implies that if $s:\O_X(-n)\ra E$ is stable then $E$
is $\tau$-semistable. The {\it class\/} of a pair $s:\O_{X}(-n)\ra
E$ is the numerical class $[E]$ in~$K(X)$. We will use $\tau'$ {\it
to denote stability of pairs}, defined using~$\O_X(1)$.
\label{ig4def3}
\end{dfn}

In \cite[Th.s~5.22 \& 5.23]{JoSo} we use results of Le Potier to
prove:

\begin{thm} If\/ $n$ is sufficiently large then the moduli functor
of stable pairs has a fine moduli scheme, a projective $\C$-scheme\/
$\M_\stp^{\al,n}(\tau'),$ with a symmetric obstruction theory.
\label{ig4thm5}
\end{thm}

\begin{dfn} In the situation above, for $\al\in K(X)$
and $n\gg 0$, define {\it stable pair invariants\/}
$PI^{\al,n}(\tau')$ in $\Z$ by
\e
\ts PI^{\al,n}(\tau')=\int_{[\M_\stp^{\al,n}(\tau')]^\vir}1,
\label{ig4eq7}
\e
where $[\M_\stp^{\al,n}(\tau')]^\vir\in A_0(\M_\stp^{\al,n}(\tau'))$
is the virtual class constructed by Behrend and Fantechi
\cite{BeFa1} using the symmetric obstruction theory from Theorem
\ref{ig4thm5}. Theorem \ref{ig3thm6}(iv) implies that the stable
pair invariants may also be written
\e
PI^{\al,n}(\tau')=\chi\bigl(\M_\stp^{\al,n}(\tau'),
\nu_{\M_\stp^{\al,n}(\tau')}\bigr).
\label{ig4eq8}
\e
\label{ig4def4}
\end{dfn}

In \cite[Cor.~5.26]{JoSo} we prove an analogue of
Theorem~\ref{ig3thm5}:

\begin{thm} $PI^{\al,n}(\tau')$ is unchanged by continuous
deformations of the underlying Calabi--Yau $3$-fold~$X$.
\label{ig4thm6}
\end{thm}

In \cite[Th.~5.27]{JoSo} we express the pair invariants
$PI^{\al,n}(\tau')$ above in terms of the generalized
Donaldson--Thomas invariants $\bar{DT}{}^\be(\tau)$ of \S\ref{ig42}.
Equation \eq{ig4eq9} is a wall-crossing formula similar to
\eq{ig4eq6}, and we prove it by change of stability condition in an
auxiliary abelian category.

\begin{thm} For\/ $\al\in C(X)$ and\/ $n\gg 0$ we have
\e
PI^{\al,n}(\tau')=\!\!\!\!\!\!\!\!\!\!
\sum_{\begin{subarray}{l} \al_1,\ldots,\al_l\in
C(X),\\ l\ge 1:\; \al_1
+\cdots+\al_l=\al,\\
\tau(\al_i)=\tau(\al),\text{ all\/ $i$}
\end{subarray} \!\!\!\!\!\!\!\!\! }
\begin{aligned}[t] \frac{(-1)^l}{l!} &\prod_{i=1}^{l}\bigl[
(-1)^{\bar\chi([\O_X(-n)]-\al_1-\cdots-\al_{i-1},\al_i)} \\
&\bar\chi\bigl([\O_{X}(-n)]\!-\!\al_1\!-\!\cdots\!-\!\al_{i-1},\al_i
\bigr) \bar{DT}{}^{\al_i}(\tau)\bigr],\!\!\!\!\!\!\!\!\!\!\!\!
\end{aligned}
\label{ig4eq9}
\e
where there are only finitely many nonzero terms in the sum.
\label{ig4thm7}
\end{thm}

Equation \eq{ig4eq9} is useful for computing invariants
$\bar{DT}{}^\al(\tau)$ in examples. By combining Theorems
\ref{ig4thm6} and \ref{ig4thm7} and using induction on the leading
coefficient of the Hilbert polynomial of $\al$, we
deduce~\cite[Cor.~5.28]{JoSo}:

\begin{cor} The generalized Donaldson--Thomas invariants
$\bar{DT}{}^\al(\tau)$ defined in\/ {\rm\S\ref{ig42}} are unchanged
under continuous deformations of the underlying Calabi--Yau\/
$3$-fold\/~$X$.
\label{ig4cor2}
\end{cor}

\subsection{Integrality properties of the $\bar{DT}{}^\al(\tau)$}
\label{ig44}

This subsection is based on ideas in Kontsevich and Soibelman
\cite[\S 2.5 \& \S 7.1]{KoSo1}. The following example is taken
from~\cite[Ex.s~6.1 \& 6.2]{JoSo}.

\begin{ex} Let $X$ be a Calabi--Yau 3-fold over $\C$ equipped
with a very ample line bundle $\O_X(1)$. Suppose $\al\in C(X)$, and
that $E\in\coh(X)$ with $[E]=\al$ is $\tau$-stable and rigid, so
that $\Ext^1(E,E)=0$. Then $mE={\buildrel {\!\ulcorner\,\text{$m$
copies }\,\urcorner\!} \over {E\op\cdots \op E}}$ for $m\ge 2$ is a
strictly $\tau$-semistable sheaf of class $m\al$, which is also
rigid. For simplicity, assume that $mE$ is the only
$\tau$-semistable sheaf of class $m\al$ for all $m\ge 1$, up to
isomorphism, so that~$\M_\ss^{m\al}(\tau)=\{[mE]\}$.

A pair $s:\O(-n)\ra mE$ may be regarded as $m$ elements
$s^1,\ldots,s^m$ of $H^0(E(n))\cong\C^{P_\al(n)}$, where $P_\al$ is
the Hilbert polynomial of $E$. Such a pair turns out to be stable if
and only if $s^1,\ldots,s^m$ are linearly independent in
$H^0(E(n))$. Two such pairs are equivalent if they are identified
under the action of $\Aut(mE)\cong\GL(m,\C)$, acting in the obvious
way on $(s^1,\ldots,s^m)$. Thus, equivalence classes of stable pairs
correspond to linear subspaces of dimension $m$ in $H^0(E(n))$, so
the moduli space $\smash{\M_\stp^{m\al,n}(\tau')}$ is isomorphic as
a $\C$-scheme to the Grassmannian $\Gr(\C^m,\C^{P_\al(n)})$. This is
smooth of dimension $m(P_\al(n)-m)$, so that
$\nu_{\M_\stp^{m\al,n}(\tau')}\equiv (-1)^{m(P_\al(n)-m)}$ by
Theorem \ref{ig3thm6}(i). Also $\Gr(\C^m,\C^{P_\al(n)})$ has Euler
characteristic the binomial coefficient
$\smash{\binom{P_\al(n)}{m}}$. Therefore \eq{ig4eq8} gives
\e
PI^{m\al,n}(\tau')=\ts (-1)^{m(P_\al(n)-m)}\binom{P_\al(n)}{m}.
\label{ig4eq10}
\e

Consider \eq{ig4eq9} with $m\al$ in place of $\al$. If
$\al_1,\ldots,\al_l$ give a nonzero term on the right hand side of
\eq{ig4eq9} then $m\al=\al_1+\cdots+\al_l$, and
$\bar{DT}{}^{\al_i}(\tau)\ne 0$, so there exists a $\tau$-semistable
$E_i$ in class $\al_i$. Thus $E_1\op\cdots\op E_l$ lies in class
$m\al$, and is $\tau$-semistable as $\tau(\al_i)=\tau(\al)$ for all
$i$. Hence $E_1\op\cdots\op E_l\cong mE$, which implies that
$E_i\cong k_iE$ for some $k_1,\ldots,k_l\ge 1$ with
$k_1+\cdots+k_l=m$, and~$\al_i=k_i\al$.

Setting $\al_i=k_i\al$, we see that $\bar\chi(\al_j,\al_i)=0$ and
$\bar\chi([\O_{X}(-n)],\al_i)=k_iP_\al(n)$, where $P_\al$ is the
Hilbert polynomial of $E$. Thus in \eq{ig4eq9} we have
$\bar\chi([\O_X(-n)]-\al_1-\cdots-\al_{i-1},\al_i)=k_iP_\al(n)$.
Combining \eq{ig4eq10}, and \eq{ig4eq9} with these substitutions,
and cancelling a factor of $(-1)^{mP_\al(n)}$ on both sides, yields
\begin{equation*}
(-1)^m\binom{P_\al(n)}{m}=
\sum_{\begin{subarray}{l}
l,k_1,\ldots,k_l\ge 1:\\
k_1+\cdots+k_l=m\end{subarray}}
\begin{aligned}[t] \frac{(-1)^l}{l!} &\prod_{i=1}^{l}
k_iP_\al(n)\bar{DT}{}^{k_i\al}(\tau).
\end{aligned}
\end{equation*}
Regarding each side as a polynomial in $P_\al(n)$ and taking the
linear term in $P_\al(n)$ we see that
\begin{equation*}
\bar{DT}{}^{m\al}(\tau)=\frac{1}{m^2} \quad\text{for all $m\ge 1$.}
\end{equation*}
\label{ig4ex1}
\end{ex}

Example \ref{ig4ex1} shows that given a rigid $\tau$-stable sheaf
$E$ in class $\al$, the sheaves $mE$ contribute $1/m^2$ to
$\bar{DT}{}^{m\al}(\tau)$ for all $m\ge 1$. We can regard this as a
kind of `multiple cover formula', analogous to the well known
Aspinwall--Morrison computation for a Calabi--Yau 3-fold $X$ that a
rigid embedded $\CP^1$ in class $\al\in H_2(X;\Z)$ contributes
$1/m^3$ to the genus zero Gromov--Witten invariant of $X$ in class
$m\al$ for all $m\ge 1$. So we can define new invariants
$\hat{DT}{}^\al(\tau)$ which subtract out these contributions from
$mE$ for~$m>1$.

\begin{dfn} Let $X$ be a projective Calabi--Yau 3-fold over $\C$,
let $\O_X(1)$ be a very ample line bundle on $X$, and let
$(\tau,T,\le)$ be a weak stability condition on $\coh(X)$ of
Gieseker or $\mu$-stability type. Then Definition \ref{ig4def2}
defines generalized Donaldson--Thomas invariants
$\bar{DT}{}^\al(\tau)\in\Q$ for~$\al\in C(X)$.

Let us define new invariants $\hat{DT}{}^\al(\tau)$ for $\al\in
C(X)$ to satisfy
\e
\bar{DT}{}^\al(\tau)=\sum_{m\ge 1,\; m\mid\al}\frac{1}{m^2}\,
\hat{DT}{}^{\al/m}(\tau).
\label{ig4eq11}
\e
By the M\"obius inversion formula, the inverse of \eq{ig4eq11} is
\e
\hat{DT}{}^\al(\tau)=\sum_{m\ge 1,\; m\mid\al}\frac{\Mo(m)}{m^2}\,
\bar{DT}{}^{\al/m}(\tau),
\label{ig4eq12}
\e
where the {\it M\"obius function\/} $\Mo:\N\ra\{-1,0,1\}$ is
$\Mo(n)=(-1)^d$ if $n=1,2,\ldots$ is square-free and has $d$ prime
factors, and $\Mo(n)=0$ otherwise.

We take \eq{ig4eq12} to be the definition of $\hat{DT}{}^\al(\tau)$,
and then reversing the argument shows that \eq{ig4eq11} holds. We
call $\hat{DT}{}^\al(\tau)$ the {\it BPS invariants\/} of $X$, as
Kontsevich and Soibelman suggest their analogous invariants
$\Om(\al)$ count BPS states.

If $\M_\ss^\al(\tau)=\M_\st^\al(\tau)$ then $\M_\ss^{\al/m}(\tau)
=\es$ for all $m\ge 2$ dividing $\al$, and so
$\hat{DT}{}^\al(\tau)=\bar{DT}{}^\al(\tau)=DT^\al(\tau)$, as in
Definition~\ref{ig4def2}.
\label{ig4def5}
\end{dfn}

We make a conjecture \cite[Conj.~6.12]{JoSo}, based
on~\cite[Conj.~6]{KoSo1}.

\begin{conj} Let\/ $X$ be a Calabi--Yau $3$-fold over $\C,$ and\/
$(\tau,T,\le)$ a weak stability condition on $\coh(X)$ of Gieseker
or $\mu$-stability type. Call\/ $(\tau,T,\le)$
\begin{bfseries}generic\end{bfseries} if for all\/ $\al,\be\in
C(X)$ with\/ $\tau(\al)=\tau(\be)$ we have\/~$\bar\chi(\al,\be)=0$.

If\/ $(\tau,T,\le)$ is generic, then $\hat{DT}{}^\al(\tau)\in\Z$ for
all\/ $\al\in C(X)$.
\label{ig4conj1}
\end{conj}

The author, Martijn Kool and Sven Meinhardt are working on a proof
of Conjecture \ref{ig4conj1}. In \cite[\S 6]{JoSo} we prove that
Conjecture \ref{ig4conj1} holds in a number of examples, and give an
example \cite[Ex.~6.8]{JoSo} in which $(\tau,T,\le)$ is not generic
and $\hat{DT}{}^\al(\tau)\notin\Z$. In \cite[Th.~7.29]{JoSo} we
prove the analogue of Conjecture \ref{ig4conj1} for invariants
counting representations of quivers without relations.

\subsection{Counting dimension 0 and 1 sheaves}
\label{ig45}

Let $X$ be a Calabi--Yau 3-fold over $\C$ with $H^1(\O_X)=0$, let
$\O_X(1)$ be a very ample line bundle on $X$, and $(\tau,G,\le)$ the
associated Gieseker stability condition on $\coh(X)$. The Chern
character gives an injective group homomorphism $\ch:K(X)\ra H^{\rm
even}(X;\Q)$. So we can regard $K(X)$ as a subgroup of $H^{\rm
even}(X;\Q)$, and write $\al\in K(X)$ as $(\al_0,\al_2,\al_4,\al_6)$
with $\al_{2j}\!\in\!H^{2j}(X;\Q)$. If $E\ra X$ is a vector bundle
with $[E]=\al$ then~$\al_0=\rank E\in\Z$.

We will consider invariants $\bar{DT}{}^\al(\tau),\hat{DT}{}^\al
(\tau)$ counting pure sheaves $E$ of dimensions 0 and 1 on $X$,
following \cite[\S 6.3--\S 6.4]{JoSo}. For sheaves $E$ of dimension
zero $\ch E=(0,0,0,d)$ where $d\ge 1$ is the length of $E$. In
\cite[\S 6.3]{JoSo} we observe that for dimension 0 sheaves the
moduli scheme $\M_\stp^{(0,0,0,d),n}(\tau')$ is independent of $n$,
and is isomorphic to the Hilbert scheme $\Hilb^dX$. Therefore
\eq{ig4eq8} gives
\begin{equation*}
PI^{(0,0,0,d),n}(\tau')=\chi\bigl(\Hilb^dX,\nu_{\Hilb^dX}\bigr),
\quad\text{for all $n\in\Z$ and $d\ge 0$.}
\end{equation*}

Values for $\chi(\Hilb^dX,\nu_{\Hilb^dX})$ were conjectured by
Maulik et al.\ \cite[Conj.~1]{MNOP1}, and proved by Behrend and
Fantechi \cite[Th.~4.12]{BeFa2} and others. These yield a generating
function for the $PI^{(0,0,0,d),n}(\tau')$:
\begin{equation*}
\ts 1+\sum_{d\ge 1} PI^{(0,0,0,d),n}(\tau')s^d= \ts\bigl[\prod_{k\ge
1}(1\!-\! s^k)^{-k}\bigr]^{\chi(X)}.
\end{equation*}
Computing using \eq{ig4eq9} then shows that
\begin{equation*}
\bar{DT}{}^{(0,0,0,d)}(\tau)=-\chi(X)\sum_{l\ge 1, \; l \mid
d}\frac{1}{l^2}.
\end{equation*}
So from \eq{ig4eq11}--\eq{ig4eq12} we deduce that
\begin{equation*}
\hat{DT}{}^{(0,0,0,d)}(\tau)=-\chi(X),\quad\text{all $d\ge 1$.}
\end{equation*}
This confirms Conjecture \ref{ig4conj1} for dimension 0 sheaves. It
is one of several examples in \cite{JoSo} in which the values of the
$PI^{\al,n}(\tau')$ are complex, the values of the
$\bar{DT}{}^\al(\tau)$ are simpler, and the values of the
$\hat{DT}{}^\al(\tau)$ are simpler still, which suggests that of the
three the invariants $\hat{DT}{}^\al(\tau)$ are the most
fundamental.

Now let $\be\in H^4(X;\Z)$ and $k\in\Z$. In \cite[\S 6.4]{JoSo} we
study invariants $\bar{DT}{}^{(0,0,\be,k)}(\tau),\hat{DT}{}^{(0,0,
\be,k)}(\tau)$ counting semistable dimension 1 sheaves, that is,
sheaves $E$ supported on curves $C$ in $X$. One expects these to be
related to curve-counting invariants like Gromov--Witten invariants,
as in the MNOP Conjecture \cite{MNOP1,MNOP2}. Here is a summary of
our results:
\begin{itemize}
\setlength{\itemsep}{0pt}
\setlength{\parsep}{0pt}
\item[(a)] $\bar{DT}{}^{(0,0,\be,k)}(\tau),\hat{DT}{}^{(0,0,
\be,k)}(\tau)$ are independent of the choice of\/~$(\tau,T,\le)$.
\item[(b)] Assume Conjecture \ref{ig4conj1} holds.
Then~$\hat{DT}{}^{(0,0,\be,k)}(\tau)\in\Z$.
\item[(c)] For any\/ $l\in \be\cup H^2(X;\Z)\subseteq\Z$ we
have $\bar{DT}{}^{(0,0,\be,k)}(\tau)=\bar{DT}{}^{(0,0,\be,k+l)}
(\tau)$ and\/ $\hat{DT}{}^{(0,0,\be,k)}(\tau)
=\hat{DT}{}^{(0,0,\be,k+l)}(\tau)$.
\item[(d)] Let $C$ be an embedded rational curve in $X$ with normal
bundle $\O(-1)\op\O(-1)$, and $\be\in H^4(X;\Z)$ be Poincar\'e dual
to $[C]\in H_2(X;\Z)$. Then sheaves supported on $C$ contribute
$1/m^2$ to $\bar{DT}{}^{(0,0,m\be,k)}(\tau)$ if $m\ge 1$ and $m\mid
k,$ and contribute $0$ to $\bar{DT}{}^{(0,0,m\be,k)}(\tau)$ if $m\ge
1$ and\/ $m\nmid k$. They contribute $1$ to
$\hat{DT}{}^{(0,0,\be,k)}(\tau),$ and $0$ to
$\hat{DT}{}^{(0,0,m\be,k)}(\tau)$ if~$m>1$.
\item[(e)] Let $C$ be a nonsingular embedded curve in $X$ of genus
$g\ge 1$, and let $\be\in H^4(X;\Z)$ be Poincar\'e dual to $[C]\in
H_2(X;\Z)$. Then sheaves supported on $C$ contribute 0 to
$\bar{DT}{}^{(0,0,m\be,k)}(\tau),\hat{DT}{}^{(0,0,m\be,k)}(\tau)$
for all $m\ge 1$ and~$k\in\Z$.
\end{itemize}
Motivated by these and by Katz \cite[Conj.~2.3]{Katz}, we
conjecture~\cite[Conj.~6.20]{JoSo}:

\begin{conj} Let\/ $X$ be a Calabi--Yau $3$-fold over\/ $\C,$ and\/
$(\tau,T,\le)$ a weak stability condition on $\coh(X)$ of Gieseker
or $\mu$-stability type. Then for $\ga\in H_2(X;\Z)$ with\/ $\be\in
H^4(X;\Z)$ Poincar\'e dual to $\ga$ and all\/ $k\in\Z$ we have
$\hat{DT}{}^{(0,0,\be,k)}(\tau)=GV_0(\ga)$. In particular,
$\hat{DT}{}^{(0,0,\be,k)}(\tau)$ is independent of\/~$k,\tau$.
\label{ig4conj2}
\end{conj}

Here $GV_0(\ga)$ is the {\it genus zero Gopakumar--Vafa invariant},
given in terms of the genus zero Gromov--Witten invariants
$GW_0(\ga)\in\Q$ of $X$ by
\begin{equation*}
GW_0(\ga)=\sum_{m\mid\ga}\frac{1}{m^3}\,GV_0(\ga/m).
\end{equation*}
A priori we have $GV_0(\ga)\in\Q$, but Gopakumar and Vafa
\cite{GoVa} conjecture that the $GV_0(\ga)$ are integers, and count
something meaningful in String Theory.

\section{Quivers with superpotentials}
\label{ig5}

We now summarize \cite[\S 7]{JoSo}, which develops an analogue of
the results of \S\ref{ig4} for representations of a quiver $Q$ with
relations $I$ coming from a superpotential $W$. In the quiver case
we have no analogue of $\bar{DT}{}^\al(\tau),\hat{DT}{}^\al(\tau),
PI^{\al,n}(\tau)$ in \S\ref{ig4} being {\it deformation-invariant},
since the proof of deformation-invariance uses the fact that the
moduli scheme $\M_\stp^{\al,n}(\tau')$ is proper (i.e.\ compact),
but the analogous moduli schemes $\smash{\M_{\stf\,Q,I}^{\bs d,\bs
e}(\mu')}$ in the quiver case need not be proper. However, all the
other important aspects of the sheaf case transfer to the quiver
case.

\subsection{Background on quivers}
\label{ig51}

Here are the basic definitions in quiver theory.

\begin{dfn} A {\it quiver\/} $Q$ is a finite directed graph.
That is, $Q$ is a quadruple $(Q_0,Q_1,h,t)$, where $Q_0$ is a finite
set of {\it vertices}, $Q_1$ is a finite set of {\it edges}, and
$h,t:Q_1\ra Q_0$ are maps giving the {\it head\/} and {\it tail\/}
of each edge.

The {\it path algebra} $\C Q$ is an associative algebra over $\C$
with basis all {\it paths of length\/} $k\ge 0$, that is, sequences
of the form
\begin{equation*}
v_0\,{\buildrel e_1\over\longra}\, v_1\ra\cdots\ra
v_{k-1}\,{\buildrel e_k\over\longra}\,v_k,
\end{equation*}
where $v_0,\ldots,v_k\in Q_0$, $e_1,\ldots,e_k\in Q_1$,
$t(a_i)=v_{i-1}$ and $h(a_i)=v_i$. Multiplication is given by
composition of paths in reverse order.

For $n\ge 0$, write $\C Q_{(n)}$ for the vector subspace of $\C Q$
with basis all paths of length $k\ge n$. It is an ideal in $\C Q$. A
{\it quiver with relations\/} $(Q,I)$ is defined to be a quiver $Q$
together with a two-sided ideal $I$ in $\C Q$ with $I\subseteq\C
Q_{(2)}$. Then $\C Q/I$ is an associative $\C$-algebra.

For $v\in Q_0$, write $i_v\in\C Q$ for the path of length 0 at $v$.
The image of $i_v$ in $\C Q/I$ is also written $i_v$. Then
\e
\text{$i_v^2=i_v,\;\> i_vi_w=0$ if $v\ne w\in Q_0,\;\>$ and
$\;\>\ts\sum_{v\in Q_0}i_v=1$ in $\C Q$ or $\C Q/I$.}
\label{ig5eq1}
\e

Write $\modCQ$ or $\modCQI$ for the abelian categories of
finite-dimensional left $\C Q$ or $\C Q/I$-modules, respectively. If
$E\in\modCQ$ or $\modCQI$ then \eq{ig5eq1} implies a decomposition
of complex vector spaces $E=\bigop_{v\in Q_0}i_v(E)$. Define the
{\it dimension vector\/} $\bdim E\in\Z^{Q_0}_{\sst\ge
0}\subset\Z^{Q_0}$ by $\bdim E:v\mapsto\dim_\C(i_vE)$. If $0\ra E\ra
F\ra G\ra 0$ is an exact sequence in $\modCQ$ or $\modCQI$ then
$\bdim F=\bdim E+\bdim G$. Hence $\bdim$ induces surjective
morphisms $\bdim:K_0(\modCQ)\ra\Z^{Q_0}$
and~$\bdim:K_0(\modCQI)\ra\Z^{Q_0}$.

Write $K(\modCQ)=K(\modCQI)=\Z^{Q_0}$, regarded as quotients of the
Grothendieck groups $K_0(\modCQ),K_0(\modCQI)$ induced by $\bdim$.
Write $C(\modCQ)=C(\modCQI)=\Z_{\sst\ge 0}^{Q_0}\sm\{0\}$, the
subsets of classes in $K(\modCQ)$, $K(\modCQI)$ of nonzero objects
in $\modCQ,\modCQI$. Here $K(\modCQ)$, $K(\modCQI)$ are our
substitutes for $K(X)=K^\num(\coh(X))$ in \S\ref{ig3}--\S\ref{ig4}.
We do not use the numerical Grothendieck groups
$K^\num(\modCQ),K^\num(\modCQI)$, as these may be zero in
interesting cases.
\label{ig5def1}
\end{dfn}

\begin{dfn} Let $Q$ be a quiver. A {\it superpotential\/} $W$ for
$Q$ over $\C$ is an element of $\C Q/[\C Q,\C Q]$. The cycles in $Q$
up to cyclic permutation form a basis for $\C Q/[\C Q,\C Q]$ over
$\C$, so we can think of $W$ as a finite $\C$-linear combination of
cycles up to cyclic permutation. We call $W$ {\it minimal\/} if all
cycles in $W$ have length at least 3. We will consider only minimal
superpotentials~$W$.

Define $I$ to be the two-sided ideal in $\C Q$ generated by $\pd_eW$
for all edges $e\in Q_1$, where if $C$ is a cycle in $Q$, we define
$\pd_eC$ to be the sum over all occurrences of the edge $e$ in $C$
of the path obtained by cyclically permuting $C$ until $e$ is in
first position, and then deleting it. Since $W$ is minimal,
$I\subseteq\C Q_{(2)}$, and $(Q,I)$ is a quiver with relations. We
allow $W\equiv 0$, so that~$I=0$.
\label{ig5def2}
\end{dfn}

Here is \cite[Th.~7.6]{JoSo}, which gives an analogue of equation
\eq{ig3eq2} for quivers with superpotentials. Now \eq{ig3eq2}
depended crucially on $X$ being a Calabi--Yau 3-fold, which implies
that $\coh(X)$ has Serre duality in dimension 3. In general the
categories $\modCQI$ coming from quivers with superpotentials do
{\it not\/} have Serre duality in dimension 3. However, as explained
in \cite[\S 7.2]{JoSo}, if $(Q,I)$ comes from a quiver with
superpotential then we can embed $\modCQI$ as the heart of a
t-structure in a 3-Calabi--Yau triangulated category $\cal T$ (which
is usually not $D^b\modCQI$), and Serre duality in dimension 3 holds
in $\cal T$. This is why quivers with superpotentials are algebraic
analogues of Calabi--Yau 3-folds, and have a version of
Donaldson--Thomas theory.

\begin{thm} Let\/ $Q=(Q_0,Q_1,h,t)$ be a quiver with relations $I$
coming from a minimal superpotential\/ $W$ on\/ $Q$ over\/ $\C$.
Define $\bar\chi:\Z^{Q_0}\t\Z^{Q_0}\ra\Z$ by
\e
\bar\chi(\bs d,\bs e)=\ts\sum_{e\in Q_1}\bigl(\bs d(h(e))\bs
e(t(e))-\bs d(t(e))\bs e(h(e))\bigr).
\label{ig5eq2}
\e
Then for any $D,E\in\modCQI$ we have
\begin{align*}
\bar\chi\bigl(\bdim D,\bdim E\bigr)=
\,&\bigl(\dim\Hom(D,E)-\dim\Ext^1(D,E)\bigr)-\\
&\bigl(\dim\Hom(E,D)-\dim\Ext^1(E,D)\bigr).
\end{align*}
\label{ig5thm1}
\end{thm}

If $Q$ is a quiver, the moduli stack $\fM_Q$ of objects $E$ in
$\modCQ$ is an Artin $\C$-stack. For $\bs d\in\Z_{\sst\ge 0}^{Q_0}$,
the open substack $\fM^{\bs d}_Q$ of $E$ with $\bdim E=\bs d$ has a
very explicit description: as a quotient $\C$-stack we have
\begin{equation*}
\fM^{\bs d}_Q\cong\ts\bigl[\prod_{e\in Q_1}\Hom(\C^{\bs
d(t(e))},\C^{\bs d(h(e))})/\prod_{v\in Q_0}\GL(\bs d(v))\bigr].
\end{equation*}
If $(Q,I)$ is a quiver with relations, the moduli stack $\fM_{Q,I}$
of objects $E$ in $\modCQI$ is a substack of $\fM_Q$, and for $\bs
d\in\Z_{\sst\ge 0}^{Q_0}$ we may write
\e
\fM^{\bs d}_{Q,I}\cong\ts\bigl[V_{Q,I}^{\bs d}/\prod_{v\in
Q_0}\GL(\bs d(v))\bigr],
\label{ig5eq3}
\e
where $V_{Q,I}^{\bs d}$ is a closed $\prod_{v\in Q_0}\GL(\bs
d(v))$-invariant $\C$-subscheme of $\prod_{e\in Q_1}\Hom\ab(\C^{\bs
d(t(e))},\C^{\bs d(h(e))})$ defined using the relations~$I$.

When $I$ comes from a superpotential $W$, we can improve the
description \eq{ig5eq3} of the moduli stacks $\fM^{\bs d}_{Q,I}$.
Define a $\prod_{v\in Q_0}\GL(\bs d(v))$-invariant polynomial
\begin{equation*}
\ts W^{\bs d}:\prod_{e\in Q_1}\Hom\bigl(\C^{\bs d(t(e))},\C^{\bs
d(h(e))}\bigr)\longra\C
\end{equation*}
as follows. Write $W$ as a finite sum $\sum_i\ga^iC^i$,where
$\ga^i\in\C$ and $C^i$ is a cycle $v_0^i\,{\buildrel
e_1^i\over\longra}\, v_1^i\ra\cdots\ra v_{k^i-1}^i\,{\buildrel
e_{k^i}^i\over\longra}\,v_{k^i}^i=v_0^i$ in $Q$. Set
\begin{equation*}
W^{\bs d}\bigl(A_e:e\in Q_1\bigr)=\ts\sum_i\ga^i
\Tr\bigl(A_{e_{k^i}^i}\ci A_{e_{k^i-1}^i}\ci\cdots \ci
A_{e_1^i}\bigr).
\end{equation*}
Then $V_{Q,I}^{\bs d}=\Crit(W^{\bs d})$ in \eq{ig5eq3}, so that
\e
\fM^{\bs d}_{Q,I}\cong\ts\bigl[\Crit(W^{\bs d})/\prod_{v\in
Q_0}\GL(\bs d(v))\bigr].
\label{ig5eq4}
\e
Equation \eq{ig5eq4} is an analogue of Theorem \ref{ig4thm1} for
categories $\modCQI$ coming from a superpotential $W$ on $Q$.

We define a class of {\it stability conditions\/} on
$\modCQI$,~\cite[Ex.~4.14]{Joyc7}.

\begin{ex} Let $(Q,I)$ be a quiver with relations. Let $c:Q_0\ra\R$
and $r:Q_0\ra(0,\iy)$ be maps. Define $\mu:C(\modCQI)\ra\R$ by
\begin{equation*}
\mu(\bs d)=\frac{\sum_{v\in Q_0}c(v)\bs d(v)}{\sum_{v\in Q_0}r(v)\bs
d(v)}\,.
\end{equation*}
Note that $\sum_{v\in Q_0}r(v)\bs d(v)>0$ as $r(v)>0$ for all $v\in
Q_0$, and $\bs d(v)\ge 0$ for all $v$ with $\bs d(v)>0$ for some
$v$. Then \cite[Ex.~4.14]{Joyc7} shows that $(\mu,\R,\le)$ is a {\it
permissible stability condition\/} on $\modCQ$, which we call {\it
slope stability}. Write $\fM_\ss^{\bs d}(\mu)$ for the open
$\C$-substack of $\mu$-semistable objects in class $\bs d$
in~$\fM^{\bs d}_{Q,I}$.

A simple case is to take $c\equiv 0$ and $r\equiv 1$, so that
$\mu\equiv 0$. Then $(0,\R,\le)$ is a trivial stability condition on
$\modCQ$ or $\modCQI$, and every nonzero object in $\modCQ$ or
$\modCQI$ is 0-semistable, so that $\fM_\ss^{\bs d}(0)=\fM^{\bs
d}_{Q,I}$.
\label{ig5ex1}
\end{ex}

\subsection{Behrend function identities, Lie
algebra morphisms, and Donaldson--Thomas type invariants}
\label{ig52}

Let $Q$ be a quiver with relations $I$ coming from a minimal
superpotential $W$ on $Q$ over $\C$. We now generalize \S\ref{ig4}
from $\coh(X)$ to $\modCQI$. The proof of Theorem \ref{ig4thm2}
depends on two things: the description of $\fM$ in terms of
$\Crit(f)$ in Theorem \ref{ig4thm1}, and equation \eq{ig3eq2}. For
$\modCQI$ equation \eq{ig5eq4} provides an analogue of Theorem
\ref{ig4thm1}, and Theorem \ref{ig5thm1} an analogue of \eq{ig3eq2}.
Thus, the proof of Theorem \ref{ig4thm2} also
yields~\cite[Th.~7.11]{JoSo}:

\begin{thm} In the situation above, with\/ $\fM_{Q,I}$ the moduli
stack of objects in a category\/ $\modCQI$ coming from a quiver $Q$
with superpotential\/ $W,$ and\/ $\bar\chi$ defined in
{\rm\eq{ig5eq2},} the Behrend function $\nu_{\fM_{Q,I}}$ of\/
$\fM_{Q,I}$ satisfies the identities \eq{ig4eq1}--\eq{ig4eq2} for
all\/~$E_1,E_2\in\modCQI$.
\label{ig5thm2}
\end{thm}

Here is the analogue of Definition~\ref{ig4def1}.

\begin{dfn} Define a Lie algebra $\ti L(Q)$ to be the $\Q$-vector
space with basis of symbols $\ti\la^{\bs d}$ for $\bs d\in\Z^{Q_0}$,
with Lie bracket
\begin{equation*}
[\ti\la^{\bs d},\ti \la^{\bs e}]=(-1)^{\bar\chi(\bs d,\bs e)}
\bar\chi(\bs d,\bs e)\ti\la^{\bs d+\bs e},
\end{equation*}
as for \eq{ig4eq3}, with $\bar\chi$ given in \eq{ig5eq2}. This makes
$\ti L(Q)$ into an infinite-dimensional Lie algebra over $\Q$.
Define $\Q$-linear maps $\ti\Psi^{\chi,\Q}_{Q,I}:\oSFai
(\fM_{Q,I},\chi,\Q)\ra\ti L(Q)$ and
$\ti\Psi_{Q,I}:\SFai(\fM_{Q,I})\ab\ra\ti L(Q)$ exactly as for
$\ti\Psi^{\chi,\Q},\ti\Psi$ in Definition~\ref{ig4def1}.
\label{ig5def3}
\end{dfn}

The proof of Theorem \ref{ig4thm3} has two ingredients: equation
\eq{ig3eq2} and Theorem \ref{ig4thm2}. Theorems \ref{ig5thm1} and
\ref{ig5thm2} are analogues of these for quivers with
superpotentials. So the proof of Theorem \ref{ig4thm3} also
yields~\cite[Th.~7.14]{JoSo}:

\begin{thm} $\ti\Psi_{Q,I}:\SFai(\fM_{Q,I})\ra\ti L(Q)$ and\/
$\ti\Psi^{\chi,\Q}_{Q,I}:\oSFai(\fM_{Q,I},\chi,\Q)\ab\ra\ti L(Q)$
are Lie algebra morphisms.
\label{ig5thm3}
\end{thm}

Here is the analogue of Definitions \ref{ig4def2} and~\ref{ig4def5}.

\begin{dfn} Let $(\mu,\R,\le)$ be a slope stability condition on
$\modCQI$ as in Example \ref{ig5ex1}. As in \S\ref{ig32} we have
elements $\bdss^{\bs d}(\mu)\in\SFa(\fM_{Q,I})$ and $\bep^{\bs
d}(\mu)\in\SFai(\fM_{Q,I})$ for all $\bs d\in C(\modCQI)$. As in
\eq{ig4eq5}, define {\it quiver generalized Donaldson--Thomas
invariants\/} $\bar{DT}{}^{\bs d}_{Q,I}(\mu)\in\Q$ for all $\bs d\in
C(\modCQI)$ by
\begin{equation*}
\ti\Psi_{Q,I}\bigl(\bep^{\bs d}(\mu)\bigr)=-\bar{DT}{}^{\bs
d}_{Q,I}(\mu)\ti\la^{\bs d}.
\end{equation*}

As in \eq{ig4eq12}, define {\it quiver BPS invariants\/}
$\hat{DT}{}^{\bs d}_{Q,I}(\mu)\in\Q$ by
\e
\hat{DT}{}^{\bs d}_{Q,I}(\mu)=\sum_{m\ge 1,\; m\mid\bs
d}\frac{\Mo(m)}{m^2}\, \bar{DT}{}^{\bs d/m}_{Q,I}(\mu),
\label{ig5eq5}
\e
where $\Mo:\N\ra\Q$ is the M\"obius function. As for \eq{ig4eq11},
the inverse of \eq{ig5eq5} is
\e
\bar{DT}{}^{\bs d}_{Q,I}(\mu)=\sum_{m\ge 1,\; m\mid\bs
d}\frac{1}{m^2}\, \hat{DT}{}^{\bs d/m}_{Q,I}(\mu).
\label{ig5eq6}
\e

If $W\equiv 0$, so that $\modCQI=\modCQ$, we write $\bar{DT}{}^{\bs
d}_Q(\mu),\hat{DT}{}^{\bs d}_Q(\mu)$ for $\smash{\bar{DT}{}^{\bs
d}_{Q,I}(\mu),\hat{DT}{}^{\bs d}_{Q,I}(\mu)}$. Note that $\mu\equiv
0$ is allowed as a slope stability condition, with every object in
$\modCQI$ 0-semistable, and is a natural choice. So we have
invariants $\bar{DT}{}^{\bs d}_{Q,I}(0), \hat{DT}{}^{\bs
d}_{Q,I}(0)$ and $\smash{\bar{DT}{}^{\bs d}_Q(0),\hat{DT}{}^{\bs
d}_Q(0)}$.
\label{ig5def4}
\end{dfn}

Here is the analogue of the integrality conjecture,
Conjecture~\ref{ig4conj1}.

\begin{conj} Call\/ $(\mu,\R,\le)$
\begin{bfseries}generic\end{bfseries} if for all\/ $\bs d,\bs e\in
C(\modCQI)$ with\/ $\mu(\bs d)=\mu(\bs e)$ we have\/~$\bar\chi(\bs
d,\bs e)=0$. If\/ $(\mu,\R,\le)$ is generic, then $\hat{DT}{}^{\bs
d}_{Q,I}(\mu)\in\Z$ for all\/~$\bs d\in C(\modCQI)$.
\label{ig5conj}
\end{conj}

In \cite[Th.~7.29]{JoSo} we prove Conjecture \ref{ig5conj} when
$W\equiv 0$, using results of Reineke \cite{Rein3}. That is, if
$\mu$ is generic we show $\hat{DT}{}^{\bs d}_Q(\mu)\in\Z$ for all
$\bs d$. In \cite[Th.~7.17]{JoSo} we prove an analogue of Theorem
\ref{ig4thm4}. It holds for arbitrary $\mu,\ti\mu$, without
requiring extra technical conditions as in Theorem~\ref{ig3thm3}.

\begin{thm} Let\/ $(\mu,\R,\le)$ and\/ $(\ti\mu,\R,\le)$ be any two
slope stability conditions on $\modCQI,$ and\/ $\bar\chi$ be as in
\eq{ig5eq2}. Then for all\/ $\bs d\in C(\modCQI)$ we have
\begin{align*}
&\bar{DT}{}^{\bs d}_{Q,I}(\ti\mu)=\\
&\!\!\!\sum_{\substack{\text{iso.}\\ \text{classes}\\
\text{of finite}\\ \text{sets $I$}}}\,\,
\sum_{\substack{\ka:I\ra C(\modCQI):\\ \sum_{i\in I}\ka(i)=\bs d}}\,\,
\sum_{\begin{subarray}{l} \text{connected,}\\
\text{simply-}\\ \text{connected}\\ \text{digraphs $\Ga$,}\\
\text{vertices $I$}\end{subarray}}
\begin{aligned}[t]
(-1)^{\md{I}-1} V(I,\Ga,\ka;\mu,\ti\mu)\cdot
\prod\nolimits_{i\in I} \bar{DT}{}^{\ka(i)}_{Q,I}(\mu)&\\
\cdot\, (-1)^{\frac{1}{2}\sum_{i,j\in
I}\md{\bar\chi(\ka(i),\ka(j))}}\cdot\! \prod\limits_{\text{edges
\smash{$\mathop{\bu}\limits^{\sst i}\ra\mathop{\bu}\limits^{\sst
j}$} in $\Ga$}\!\!\!\!\!\!\!\!\!\!\!\!\!\!\!\!\!\!\!\!\!\!\!
\!\!\!\!\!\!\!\!} \bar\chi(\ka(i),\ka(j))&,
\end{aligned}
\end{align*}
with only finitely many nonzero terms.
\label{ig5thm4}
\end{thm}

\subsection{Pair invariants for quivers}
\label{ig53}

We now discuss analogues for quivers of the moduli spaces of stable
pairs $\M_\stp^{\al,n}(\tau')$ and stable pair invariants
$PI^{\al,n}(\tau')$ in \S\ref{ig43}, and the identity \eq{ig4eq9} in
Theorem \ref{ig4thm7} relating $PI^{\al,n}(\tau')$ and the
$\bar{DT}{}^\be(\tau)$.

\begin{dfn} Let $Q$ be a quiver with relations $I$ coming from a
superpotential $W$ on $Q$ over $\C$. Suppose $(\mu,\R,\le)$ is a
slope stability condition on $\modCQI$, as in Example~\ref{ig5ex1}.

Let $\bs d,\bs e\in\Z^{Q_0}_{\sst\ge 0}$ be dimension vectors. A
{\it framed representation\/ $(E,\si)$ of\/ $(Q,I)$ of type\/} $(\bs
d,\bs e)$ consists of a representation $E$ of $\modCQ/I$ with $\bdim
E=\bs d$, together with linear maps $\si_v:\C^{\bs e(v)}\ra i_v(E)$
for all $v\in Q_0$. We call a framed representation $(E,\si)$ {\it
stable\/} if
\begin{itemize}
\setlength{\itemsep}{0pt}
\setlength{\parsep}{0pt}
\item[(i)] $\mu([E'])\le\mu([E])$ for all
subobjects $0\ne E'\subset E$ in $\modCQI$; and
\item[(ii)] If also $\si$ factors through $E'$, that is,
$\si_v(\C^{e(v)})\subseteq i_v(E')\subseteq i_v(E)$ for all $v\in
Q_0$, then~$\mu([E'])<\mu([E])$.
\end{itemize}
We will use $\mu'$ {\it to denote stability of framed
representations}, defined using~$\mu$.
\label{ig5def5}
\end{dfn}

Following Engel and Reineke \cite[\S 3]{EnRe} or Szendr\H oi
\cite[\S 1.2]{Szen}, we can in a standard way define moduli problems
for all framed representations, and for stable framed
representations. The moduli space of all framed representations of
type $(\bs d,\bs e)$ is an Artin $\C$-stack $\fM_{\fr\,Q,I}^{\bs
d,\bs e}$ with an explicit description similar to \eq{ig5eq4}, and
the moduli space of stable framed representations of type $(\bs
d,\bs e)$ is a fine moduli $\C$-scheme $\smash{\M_{\stf\,Q,I}^{\bs
d,\bs e}(\mu')}$, an open $\C$-substack of~$\fM_{\fr\,Q,I}^{\bs
d,\bs e}$.

We can now define our analogues of invariants $PI^{\al,n}(\tau')$
for quivers.

\begin{dfn} In the situation above, define
\e
NDT_{Q,I}^{\bs d,\bs e}(\mu')=\chi\bigl(\M_{\stf\,Q,I}^{\bs d,\bs
e}(\mu'),\nu_{\M_{\stf\,Q,I}^{\bs d,\bs e}(\mu')}\bigr).
\label{ig5eq7}
\e
When $W\equiv 0$, so that $\modCQI=\modCQ$, we also write
$NDT_{Q}^{\bs d,\bs e}(\mu')=NDT_{Q,I}^{\bs d,\bs e}(\mu')$.
Following Szendr\H oi \cite{Szen} we call $NDT_{Q,I}^{\bs d,\bs
e}(\mu'),NDT_{Q}^{\bs d,\bs e}(\mu')$ {\it noncommutative
Donaldson--Thomas invariants}.
\label{ig5def6}
\end{dfn}

Here \eq{ig5eq7} is the analogue of \eq{ig4eq8} in the sheaf case.
We have no analogue of \eq{ig4eq7}, since in general
$\M_{\stf\,Q,I}^{\bs d,\bs e}(\mu')$ is not proper, and so does not
have a fundamental class. These quiver analogues of
$\M_\stp^{\al,n}(\tau'),PI^{\al,n}(\tau')$ are not new, similar
things have been studied in quiver theory by Nakajima, Reineke,
Szendr\H oi and other authors for some years
\cite{EnRe,Naga,NaNa,Naka,MoRe,Rein1,Rein2,Szen}. Here
\cite[Th.~7.23]{JoSo} is the analogue of Theorem \ref{ig4thm7} for
quivers.

\begin{thm} Suppose\/ $Q$ is a quiver with relations\/ $I$ coming
from a minimal superpotential\/ $W$ on $Q$ over $\C$. Let\/
$(\mu,\R,\le)$ be a slope stability condition on\/ $\modCQI,$ as in
Example {\rm\ref{ig5ex1},} and\/ $\bar\chi$ be as in \eq{ig5eq2}.
Then for all\/ $\bs d,\bs e$ in\/ $C(\modCQI)=\Z_{\sst\ge
0}^{Q_0}\sm\{0\}\subset \Z^{Q_0},$ we have
\e
NDT^{\bs d,\bs e}_{Q,I}(\mu')=\!\!\!\!\!\!\!\!\!\!\!\!\!\!\!
\sum_{\begin{subarray}{l} \bs d_1,\ldots,\bs d_l\in
C(\modCQI),\\ l\ge 1:\; \bs d_1
+\cdots+\bs d_l=\bs d,\\
\mu(\bs d_i)=\mu(\bs d),\text{ all\/ $i$}
\end{subarray} \!\!\!\!\!\!\!\!\! }
\begin{aligned}[t] \frac{(-1)^l}{l!} &\prod_{i=1}^{l}\bigl[
(-1)^{\bs e\cdot\bs d_i-\bar\chi(\bs d_1+\cdots+\bs
d_{i-1},\bs d_i)} \\
&\bigl(\bs e\cdot\bs d_i-\bar\chi(\bs d_1\!+\!\cdots\!+\!\bs
d_{i-1},\bs d_i)\bigr)\bar{DT}{}^{\bs
d_i}_{Q,I}(\mu)\bigr],\!\!\!\!\!\!\!\!\!\!\!\!\!\!
\end{aligned}
\label{ig5eq8}
\e
with\/ $\bs e\cdot\bs d_i=\sum_{v\in Q_0}\bs e(v)\bs d_i(v),$ and\/
$\bar{DT}{}^{\bs d_i}_{Q,I}(\mu),NDT^{\bs d,\bs e}_{Q,I}(\mu')$ as
in Definitions\/ {\rm\ref{ig5def4}, \rm\ref{ig5def6}}. When $W\equiv
0,$ the same equation holds for $\smash{NDT^{\bs d,\bs
e}_Q(\mu'),\bar{DT}{}^{\bs d}_Q(\mu)}$.
\label{ig5thm5}
\end{thm}

For Donaldson--Thomas invariants in \S\ref{ig4}, we regarded the
invariants $\bar{DT}{}^\al(\tau)$, $\hat{DT}{}^\al(\tau)$ as our
primary objects of study, and the pair invariants
$PI^{\al,n}(\tau')$ as secondary, not of that much interest in
themselves. In contrast, in the quiver literature to date the
invariants $\bar{DT}{}^{\bs d}_{Q,I}(\mu),\bar{DT}{}^{\bs d}_Q(\mu)$
and $\hat{DT}{}^{\bs d}_{Q,I}(\mu),\hat{DT}{}^{\bs d}_Q(\mu)$ have
not been seriously considered even in the stable$\,=\,$semistable
case, and the analogues $\smash{NDT_{Q,I}^{\bs d,\bs e}(\mu'),
NDT_{Q}^{\bs d,\bs e}(\mu')}$ of pair invariants $PI^{\al,n}(\tau')$
have been the central object of study.

We argue that the invariants $\smash{\bar{DT}{}^{\bs
d}_{Q,I}(\mu),\ldots,\hat{DT}{}^{\bs d}_Q(\mu)}$ should actually be
regarded as more fundamental and more interesting than the
$NDT_{\smash{Q,I}}^{\bs d,\bs e}(\mu')$. By \eq{ig5eq8} the
$NDT_{\smash{Q,I}}^{\bs d,\bs e}(\mu')$ can be written in terms of
the $\bar{DT}{}^{\bs d}_{\smash{Q,I}}(\mu)$, and hence by
\eq{ig5eq6} in terms of the $\hat{DT}{}^{\bs d}_{\smash{Q,I}}(\mu)$,
so the pair invariants contain no more information. The
$\bar{DT}{}^{\bs d}_{\smash{Q,I}}(\mu)$ are simpler than the
$NDT_{\smash{Q,I}}^{\bs d,\bs e}(\mu')$ as they depend only on $\bs
d$ rather than on $\bs d,\bs e$. In examples in \cite[\S 7.5--\S
7.6]{JoSo} we find that the values of the $\bar{DT}{}^{\bs
d}_{\smash{Q,I}}(\mu)$ and especially of the $\hat{DT}{}^{\bs
d}_{\smash{Q,I}}(\mu)$ may be much simpler and more illuminating
than the values of the $NDT_{\smash{Q,I}}^{\bs d,\bs e}(\mu')$, as
in \eq{ig5eq9}--\eq{ig5eq11} below.

Here is an example taken from~\cite[\S 7.5.2]{JoSo}.

\begin{ex} Following Szendr\H oi \cite[\S 2.1]{Szen}, let
$Q=(Q_0,Q_1,h,t)$ have two vertices $Q_0=\{v_0,v_1\}$ and edges
$e_1,e_2:v_0\ra v_1$ and $f_1,f_2:v_1\ra v_0$, as below:
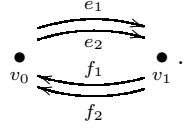
\begin{equation*}
\xymatrix@C=40pt{ \mathop{\bu} \limits_{v_0} \ar@/^/@<2ex>[r]_{e_2}
\ar@/^/@<3ex>[r]^{e_1} & \mathop{\bu}\limits_{v_1}.
\ar@/^/@<0ex>[l]_{f_1} \ar@/^/@<1ex>[l]^{f_2} }
\end{equation*}
Define a superpotential $W$ on $Q$ by $W=e_1f_1e_2f_2-e_1f_2e_2f_1$,
and let $I$ be the associated relations. Then $\modCQI$ is a
3-Calabi--Yau abelian category. Theorem \ref{ig5thm1} shows that the
Euler form $\bar\chi$ on $\modCQI$ is zero.

Write elements $\bs d$ of $C(\modCQI)$ as $(d_0,d_1)$ where $d_j=\bs
d(v_j)$. Szendr\H oi \cite[Th.~2.7.1]{Szen} computed the
noncommutative Donaldson--Thomas invariants
$NDT^{\smash{(d_0,d_1),(1,0)}}_{Q,I}(0')$ for $\modCQI$ as
combinatorial sums, and using work of Young \cite{Youn} wrote their
generating function as a product \cite[Th.~2.7.2]{Szen}, giving
\e
\begin{split}
1&+\sum_{(0,0)\ne(d_0,d_1)\in\N^2}NDT^{(d_0,d_1),(1,0)}_{Q,I}(0')
q_0^{d_0}q_1^{d_1}\\
&=\prod_{k\ge 1}\bigl(1-(-q_0q_1)^k)\bigr)^{-2k}\bigl(1-(-q_0)^k
q_1^{k-1}\bigr)^k\bigl(1-(-q_0)^kq_1^{k+1}\bigr)^k.
\end{split}
\label{ig5eq9}
\e
Computing using \eq{ig5eq8} and \eq{ig5eq9} shows that
\e
\bar{DT}{}^{(d_0,d_1)}_{Q,I}(0)=\begin{cases} \displaystyle
-2\sum_{l\ge 1, \; l \mid d}\frac{1}{l^2},\!\! &
d_0=d_1=d\ge 1, \\
\displaystyle \frac{1}{l^2}, & d_0=kl, \; d_1=(k-1)l,\;
k,l\ge 1, \\[7pt]
\displaystyle \frac{1}{l^2}, & d_0=kl, \; d_1=(k+1)l,\;
k\ge 0,\; l\ge 1,\!\! \\
0, & \text{otherwise.}
\end{cases}
\label{ig5eq10}
\e
Combining \eq{ig5eq6} and \eq{ig5eq10} we see that
\e
\hat{DT}{}^{(d_0,d_1)}_{Q,I}(0)=\begin{cases} -2, &
(d_0,d_1)=(k,k),\; k\ge 1,\\
\phantom{-}1, & (d_0,d_1)=(k,k-1),\; k\ge 1,\\
\phantom{-}1, &  (d_0,d_1)=(k-1,k),\; k\ge 1,\\
\phantom{-}0, & \text{otherwise.}
\end{cases}
\label{ig5eq11}
\e
Note that the values of the
$\smash{\hat{DT}{}^{(d_0,d_1)}_{Q,I}(0)}$ in \eq{ig5eq11} lie in
$\Z$, as in Conjecture \ref{ig5conj}, and are far simpler than those
of the $NDT^{(d_0,d_1),(1,0)}_{Q,I}(0')$ in~\eq{ig5eq9}.

This example is connected to Donaldson--Thomas theory for
(noncompact) Calabi--Yau 3-folds as follows. We have equivalences of
derived categories
\e
D^b(\modCQI)\sim D^b(\coh_\cs(X))\sim D^b(\coh_\cs(X_+)),
\label{ig5eq12}
\e
where $\pi:X\ra Y$ and $\pi_+:X_+\ra Y$ are the two crepant
resolutions of the conifold $Y=\bigl\{(z_1,z_2,z_3,z_4)\in
\C^4:z_1^2+\cdots+z_4^2=0\bigr\}$, and $X,X_+$ are related by a
flop. Here $X,X_+$ are regarded as `commutative' crepant resolutions
of $Y$, and $\modCQI$ as a `noncommutative' resolution of $Y$, in
the sense that $\modCQI$ can be regarded as the coherent sheaves on
the `noncommutative scheme' $\Spec(\C Q/I)$ constructed from the
noncommutative $\C$-algebra~$\C Q/I$.

One idea in \cite{Szen} is that counting invariants $NDT_{Q,I}^{\bs
d,\bs e}(\mu')$ for $\modCQI$ should be related to Donaldson--Thomas
type invariants counting sheaves on $X,X_+$ by some kind of
wall-crossing formula under change of stability condition in the
derived categories, using the equivalences \eq{ig5eq12}. This
picture has been worked out further by Nagao and Nakajima
\cite{NaNa,Naga}. In \cite[\S 7.5.2]{JoSo} we show that in this case
the situation for invariants $\bar{DT}{}^{\bs d}_{Q,I}(\mu),
\hat{DT}{}^{\bs d}_{Q,I}(\mu)$ is actually much simpler, because
they are unchanged by wall-crossing as $\bar\chi\equiv 0$, so we can
identify the invariants $\bar{DT}{}^{\bs d}_{Q,I}(\mu),
\hat{DT}{}^{\bs d}_{Q,I}(\mu)$ in \eq{ig5eq10}--\eq{ig5eq11}
directly with Donaldson--Thomas invariants for $X$ and~$X_+$.
\label{ig5ex2}
\end{ex}

\medskip

\noindent{\small\sc The Mathematical Institute, 24-29 St. Giles,
Oxford, OX1 3LB, U.K.}

\noindent{\small\sc E-mail: \tt joyce@maths.ox.ac.uk}


\begin{thebibliography}{99}

\bibitem{Behr} K. Behrend, {\it Donaldson--Thomas type invariants via
microlocal geometry}, to appear in Annals of Mathematics.
math.AG/0507523, 2005.

\bibitem{BeFa1} K. Behrend and B. Fantechi, {\it The intrinsic normal
cone}, Invent. Math. 128 (1997), 45--88.

\bibitem{BeFa2} K. Behrend and B. Fantechi, {\it Symmetric
obstruction theories and Hilbert schemes of points on threefolds},
Algebra and Number Theory 2 (2008), 313--345. math.AG/0512556.

\bibitem{DoTh} S.K. Donaldson and R.P. Thomas, {\it Gauge Theory in
Higher Dimensions}, Chapter 3 in S.A. Huggett, L.J. Mason, K.P. Tod,
S.T. Tsou and N.M.J. Woodhouse, editors, {\it The Geometric
Universe}, Oxford University Press, Oxford, 1998.

\bibitem{EnRe} J. Engel and M. Reineke, {\it Smooth models of
quiver moduli}, Math. Z. 262 (2009), 817--848. arXiv:0706.4306.

\bibitem{Gome} T.L. G\'omez, {\it Algebraic stacks}, Proc. Indian Acad.
Sci. Math. Sci. 111 (2001), 1--31. math.AG/9911199.

\bibitem{GoVa} R. Gopakumar and C. Vafa, {\it M-theory and
topological strings. II}, hep-th/9812127, 1998.

\bibitem{HuLe} D. Huybrechts and M. Lehn, {\it The geometry of
moduli spaces of sheaves}, Aspects of Math. E31, Vieweg,
Braunschweig/Wiesbaden, 1997.

\bibitem{Joyc1} D. Joyce, {\it Constructible functions on Artin
stacks}, J. London Math. Soc. 74 (2006), 583--606. math.AG/0403305.

\bibitem{Joyc2} D. Joyce, {\it Motivic invariants of Artin stacks
and `stack functions'}, Quart. J. Math. 58 (2007), 345--392.
math.AG/0509722.

\bibitem{Joyc3} D. Joyce, {\it Configurations in abelian
categories. I. Basic properties and moduli stacks}, Adv. Math. 203
(2006), 194--255. math.AG/0312190.

\bibitem{Joyc4} D. Joyce, {\it Configurations in abelian
categories. II. Ringel--Hall algebras}, Adv. Math. 210 (2007),
635--706. math.AG/0503029.

\bibitem{Joyc5} D. Joyce, {\it Configurations in abelian
categories. III. Stability conditions and identities}, Adv. Math.
215 (2007), 153--219. math.AG/0410267.

\bibitem{Joyc6} D. Joyce, {\it Configurations in abelian categories.
IV. Invariants and changing stability conditions}, Adv. Math. 217
(2008), 125--204. math.AG/0410268.

\bibitem{Joyc7} D. Joyce, {\it Holomorphic generating functions for
invariants counting coherent sheaves on Calabi--Yau $3$-folds},
Geometry and Topology 11 (2007), 667--725. hep-th/0607039.

\bibitem{JoSo} D. Joyce and Y. Song, {\it A theory of generalized
Donaldson--Thomas invariants}, arXiv:0810.5645v5, May 2010. 211
pages.

\bibitem{Katz} S. Katz, {\it Genus zero Gopakumar--Vafa invariants
of contractible curves}, J. Diff. Geom. 79 (2008), 185-195.
math.AG/0601193.

\bibitem{KoSo1} M. Kontsevich and Y. Soibelman, {\it Stability
structures, motivic Donaldson--Thomas invariants and cluster
transformations}, arXiv:0811.2435, 2008.

\bibitem{KoSo2} M. Kontsevich and Y. Soibelman, {\it Motivic
Donaldson--Thomas invariants: summary of results}, arXiv:0910.4315,
2009.

\bibitem{LaMo} G. Laumon and L. Moret-Bailly, {\it Champs
alg\'ebriques}, Ergeb. der Math. und ihrer Grenzgebiete 39,
Springer-Verlag, Berlin, 2000.

\bibitem{Miya} K. Miyajima, {\it Kuranishi family of vector bundles
and algebraic description of Einstein--Hermitian connections}, Publ.
RIMS, Kyoto Univ. 25 (1989), 301--320.

\bibitem{MNOP1} D. Maulik, N. Nekrasov, A. Okounkov, and R.
Pandharipande, {\it Gromov-Witten theory and Donaldson--Thomas
theory. I}, Compos. Math. 142 (2006), 1263--1285. math.AG/0312059.

\bibitem{MNOP2} D. Maulik, N. Nekrasov, A. Okounkov, and R.
Pandharipande, {\it Gromov-Witten theory and Donaldson--Thomas
theory. II}, Compos. Math. 142 (2006), 1286--1304. math.AG/0406092.

\bibitem{MoRe} S. Mozgovoy and M. Reineke, {\it On the
noncommutative Donaldson-Thomas invariants arising from brane
tilings}, arXiv:0809.0117, 2008.

\bibitem{Naga} K. Nagao, {\it Derived categories of small toric
Calabi--Yau $3$-folds and curve counting invariants},
arXiv:0809.2994, 2008.

\bibitem{NaNa} K. Nagao and H. Nakajima, {\it Counting invariant of
perverse coherent sheaves and its wall-crossing}, arXiv:0809.2992,
2008.

\bibitem{Naka} H. Nakajima, {\it Varieties associated with quivers},
pages 139--157 in R. Bautista et al., editors, {\it Representation
theory of algebras and related topics}, C.M.S. Conf. Proc. 19,
A.M.S., Providence, RI, 1996.

\bibitem{PaTh} R. Pandharipande and R.P. Thomas, {\it Curve
counting via stable pairs in the derived category}, arXiv:0707.2348,
2007.

\bibitem{Rein1} M. Reineke, {\it Cohomology of noncommutative
Hilbert schemes}, Algebr. Represent. Theory 8 (2005), 541--561.
math.AG/0306185.

\bibitem{Rein2} M. Reineke, {\it Framed quiver moduli, cohomology,
and quantum groups}, J. Algebra 320 (2008), 94--115.
math.AG/0411101.

\bibitem{Rein3} M. Reineke, {\it Cohomology of quiver moduli,
functional equations, and integrality of Donaldson--Thomas type
invariants}, arXiv:0903.0261, 2009.

\bibitem{Szen} B. Szendr\H oi, {\it Non-commutative Donaldson--Thomas
theory and the conifold}, Geom. Topol. 12 (2008), 1171--1202.
arXiv:0705.3419.

\bibitem{Thom} R.P. Thomas, {\it A holomorphic Casson invariant for
Calabi--Yau $3$-folds, and bundles on $K3$ fibrations}, J. Diff.
Geom. 54 (2000), 367--438. \hfil\break math.AG/9806111.

\bibitem{Youn} B. Young, {\it Computing a pyramid partition
generating function with dimer shuffling}, J. Combin. Theory Ser. A
116 (2009), 334--350. arXiv:0709.3079.

\end{thebibliography}
\end{document}